\documentclass[12pt,a4paper]{amsart}
\makeatletter
\renewcommand\normalsize{%
    \@setfontsize\normalsize{11.7}{14pt plus .3pt minus .3pt}%
    \abovedisplayskip 10\p@ \@plus4\p@ \@minus4\p@
    \abovedisplayshortskip 6\p@ \@plus2\p@
    \belowdisplayshortskip 6\p@ \@plus2\p@
    \belowdisplayskip \abovedisplayskip}
\renewcommand\small{%
    \@setfontsize\small{9.5}{12\p@ plus .2\p@ minus .2\p@}%
    \abovedisplayskip 8.5\p@ \@plus4\p@ \@minus1\p@
    \belowdisplayskip \abovedisplayskip
    \abovedisplayshortskip \abovedisplayskip
    \belowdisplayshortskip \abovedisplayskip}
\renewcommand\footnotesize{%
    \@setfontsize\footnotesize{8.5}{9.25\p@ plus .1pt minus .1pt}
    \abovedisplayskip 6\p@ \@plus4\p@ \@minus1\p@
    \belowdisplayskip \abovedisplayskip
    \abovedisplayshortskip \abovedisplayskip
    \belowdisplayshortskip \abovedisplayskip}
\ifdefined\pdfpagewidth
\else
\fi
\calclayout
\makeatother

\usepackage[utf8]{inputenc}
\usepackage{mathtools}
\usepackage{amssymb}
\usepackage{amsfonts}
\usepackage{amsthm}
\usepackage{bbm}
\usepackage{fullpage}
\usepackage{color}
\usepackage{hyperref}
\usepackage{stmaryrd}
\usepackage{enumitem}
\usepackage{accents}
\mathtoolsset{showonlyrefs}

\newcommand{\ind}{\mathbbm 1}

\newtheorem{theorem}{Theorem}[section]
\newtheorem{lemma}[theorem]{Lemma}
\newtheorem{proposition}[theorem]{Proposition}
\newtheorem{corollary}[theorem]{Corollary}

\newtheorem{thmx}{Theorem}

\newtheorem{remark}[theorem]{Remark}

\theoremstyle{definition}

\newcommand{\dd}{\text{d}}

\newcommand{\f}{\mathfrak{f}}

\newcounter{constants}

\renewcommand{\tilde}{\widetilde}

\newcommand{\bbE}{{\ensuremath{\mathbb E}} }

\newcommand{\bbN}{{\ensuremath{\mathbb N}} }

\newcommand{\bbP}{{\ensuremath{\mathbb P}} }

\newcommand{\bbR}{{\ensuremath{\mathbb R}} }

\newcommand{\bbZ}{{\ensuremath{\mathbb Z}} }

\renewcommand{\epsilon}{\varepsilon}

\newcommand{\gep}{\varepsilon}       

\newcommand{\go}{\omega}

\newcommand{\gl}{\lambda}

\newcommand{\cA}{{\ensuremath{\mathcal A}} }
\newcommand{\cB}{{\ensuremath{\mathcal B}} }
\newcommand{\cF}{{\ensuremath{\mathcal F}} }
\newcommand{\cJ}{{\ensuremath{\mathcal J}} }

\newcommand{\cN}{{\ensuremath{\mathcal N}} }

\newcommand{\cZ}{{\ensuremath{\mathcal Z}} }

\newcommand{\cK}{{\ensuremath{\mathcal K}} }

\newcommand{\bP}{{\ensuremath{\mathbf P}} }

\newcommand{\lint}{\llbracket}
\newcommand{\rint}{\rrbracket}

\definecolor{darkred}{rgb}{0.7,0.1,0.1}
\renewcommand{\hat}{\widehat}

\renewcommand{\complement}{\mathsf{c}}
\newcommand{\tf}{\textsc{f}}

\title{The localization transition for the directed polymer in a random environment is smooth}

\author{Hubert Lacoin}
\address{IMPA, Estrada Dona Castorina 110,
Rio de Janeiro RJ-22460-320- Brasil}
\email{lacoin@impa.br}

\begin{document}

\begin{abstract}
When $d\ge 3$, the directed polymer  in a  random environment on $\bbZ^d$ displays a phase transition from a diffusive phase, known as \textit{weak disorder}, to a localized phase, referred to as \textit{strong disorder}.
This transition is encoded by the behavior of the free energy of the model, defined by 
$$\f(\beta):=\lim_{N\to \infty} (1/n)\log W^{\beta}_n$$ where $W^{\beta}_n$ is the normalized partition function for the directed polymer of length $n$. More precisely  weak disorder corresponds to $\f(\beta)=0$ and  strong disorder to $\f(\beta)<0$. Monotonicity and continuity of $\f$ imply that there exists $\beta_c\in [0,\infty]$ such that weak disorder is equivalent to $\beta\in [0,\beta_c]$. Furthermore $\beta_c>0$ if and only if $d\ge 3$.
We prove that this transition is infinitely smooth in the sense that  $\f$ grows slower than any power function at the vicinity of $\beta_c$, that is
$$ \lim_{\beta \downarrow \beta_c }\frac{\log |\f(\beta)|}{\log (\beta-\beta_c)}=\infty.$$
\end{abstract}

\maketitle

\section{Introduction}

In statistical mechanics, the notion of \textit{phase transition} refers to an abrupt change of behavior of a physical system when a parameter passes a critical threshold. Mathematically speaking, since the Gibbs weights are analytic functions, there is no phase transition in a finite state space and the phenomenon appears  when considering infinite volume limits. The signature of a phase transition is the loss of analytic behavior of the \textit{free energy density} -- the  $\log$ of the  partition function divided by the volume -- when the size of the system goes to infinity. The understanding of the behavior of the free energy around critical points -- in particar the identification of \textit{critical exponents} -- is one of the central questions in the study  phase transitions.

\medskip

The directed polymer in a random environment on $\bbZ^d$ is a disordered model which has been continuously studied by  physicists and mathematicians for four decades. In spite of this undiscontinued interest, no mathematical results on the critical behavior of the free energy  has yet  been proved for $d\ge 3$ (when $d\le 2$ the model does not display a phase transition). In the present paper we fill this gap by establishing rigorously that the critical exponent associated to the free energy is infinite.

\medskip

The directed polymer in a random environment was first introduced (with transversal dimension $d=1$) as a toy model to describe interfaces for the low temperature Ising model with edge disorder \cite{HH85}.
The interest in the model quickly extended to higher dimensions both in the theoretical physics and mathematics communities  (see \cite{B89,DE92,IS88,KZ87} and references therein). We refer the reader to the monograph \cite{Com17} and to the recent survey paper \cite{Zyg24} for a complete introduction to the subject.

\medskip

It was proved in \cite{B89,IS88} that when $d$ is larger than $3$, the system displays a diffusive behavior at high temperature, hinting at the existence of a phase transition. The existence of a low temperature phase -- with an explicit lower bound for $\beta_c$ -- is established in \cite{CD89} while  the existence of the free energy is proved in \cite{CH02,CSY03}. A monotonicity result which garanties the uniqueness of the critical temperature is presented in  \cite{CY06}.

\medskip

To this day,  very little has been achieved, even at the heuristic level, in the  understanding  of the critical behavior of the free energy curve. The few  numerical and theoretical predictions concerning the critical exponent published in the physics literature are restricted to the case $d=3$ and illustrate that no consensus has been reached on the topic \cite{DG90,KBM91, MG2, MG1}.

\subsection{The directed polymer in a random environment and its free energy}

To push the discussion further,  we  need to introduce the model formally. 
For the sake of this introduction, we condense existing results about the localization phase transition in three theorems labeled \ref{leA},\ref{leB} and \ref{leC}). 
 We let $(X_n)_{n\ge 0}$  denote a symmetric simple random walk on $\bbZ^d$ starting from the origin and $P$ its distribution. We have $P(X_0=0)=1$ and the increments $(X_n-X_{n-1})_{n\ge 1}$ are i.i.d.\ with distribution $P(X_1=x)= \ind_{\{|x|=1\}}/(2d)$ where $|\cdot|$ denote the $\ell_1$ norm on $\bbZ^d$.
An environment $\go= (\go_{n,x})_{n\ge 1,x\in \bbZ^d}$ is a collection of real numbers indexed by $\bbN\times \bbZ^d$. Given $\go$, $\beta\ge 0$ and $n\ge 0$, the directed polymer of length $n$, with environment $\go$ and inverse temperature $\beta$ is the probability on the set of random walk trajectories whose density with respect to $P$ is given by
\begin{equation}
 \frac{\dd P^{\beta,\go}_{n}}{\dd P}(X)=\frac{1}{Z^{\beta,\go}_n} \exp\left(\beta \sum_{j=1}^n \go_{j,X_j} \right) \quad   \text{ with } \quad   Z^{\beta,\go}_n:= E\left[e^{\beta \sum_{j=1}^n \go_{j,X_j}}   \right].
\end{equation}
We consider the case where  environment $\go$  is obtained  by sampling a  collection of i.i.d.\ random variables. We let $\bbP$ denote the associated distribution. We assume that these variables admit exponential moments of all order (a standard assumption in the directed polymer literature)
\begin{equation}\label{allmoments}
 \forall \beta\in \bbR, \quad \gl(\beta):=\log \bbE\left[ e^{\beta \go_{1,0}}\right]<\infty.
\end{equation}
The phase transition for the directed polymer in a random environment is encoded by the asymptotic behavior in $n$ of the sequence of  partition functions $(Z^{\beta,\go}_n)_{n\ge 1}$. We consider the normalized partition function $W^{\beta}_n$ (the dependence in $\go$ is omitted in the notation for better readability) defined by
\begin{equation}
W^{\beta}_n:=Z^{\beta,\go}_n / \bbE[Z^{\beta,\go}_n]= E \left[ e^{\sum_{j=1}^n(\beta \go_{j,X_j}-\gl(\beta))}\right].
\end{equation}
If $d\ge 3$, the directed polymer undergoes a transition when $\beta$ grows, from a regime where $W^{\beta}_n$ converges to a nontrivial random limit $W^{\beta}_{\infty}>0$, to one where $W^{\beta}_n$ decays exponentially to zero. We summarize what is known on this topic in the next two theorems.
The first one establishes the existence and monotonicity of the free energy.

 \begin{thmx}[\cite{B89,CH02,CSY03,CV06,CY06,IS88,L10}]\label{leA}
The following limits are  well defined and coincide almost-surely
\begin{equation}\label{freenexi}
\f(\beta):=\lim_{n\to \infty}  \frac{1}{n} \bbE\left[ \log W^{\beta}_n\right]= \lim_{n\to \infty}  \frac{1}{n} \log W^{\beta}_n.
 \end{equation}
The function $\f(\beta)$ is referred to as the free energy of the directed polymer, it is continous and nonincreasing in $\beta$. Thus
there exists $\beta_c\in [0,\infty]$ such that
$$ \beta>\beta_c \quad \Leftrightarrow  \quad \f(\beta)<0.$$
Furthermore the following holds:
\begin{itemize}
 \item [(i)] The function $\f(\beta)+\gl(\beta)=\lim_{n\to \infty} (1 /n)\log Z^{\beta,\go}_n$ is convex in $\beta$.
 \item [(ii)] If $d=1,2$ we have $\beta_c=0$, while if $d\ge 3$ we have $\beta_c>0$.
\item [(iii)] When $d\ge 3$ we have $\beta_c\in (0,\infty)$  if $\bbP\left( \go_{1,0}=\mathrm{ess\, sup}  \ \go_{1,0}\right)<\frac{1}{2d}$.
\end{itemize}
\end{thmx}
The existence and coincidence of the limits in \eqref{freenexi} was proved in \cite{CH02,CSY03}, while the monotonicity in $\beta$ was established in \cite{CY06}. The identity $\beta_c=0$ was proved in \cite{CV06} for $d=1$ and in \cite{L10} for $d=2$.  For $d\ge 3$, it was proved that  $\beta_c>0$ in \cite{B89,IS88}.
The criterion to have $\beta_c<\infty$ in item $(iii)$ is established in \cite[Theorem 2.3, item (a)]{CSY03}. The value $1/2d$ is not optimal but the purpose of $(iii)$ is to illustrate the fact that, while $\beta_c=\infty$ occurs for some distribution of $\go$, it is the exception rather than the rule.
The second theorem establishes the sharpness of the phase transition, showing that $W^{\beta}_n$ admits a non trivial limit when $\beta\in [0,\beta_c]$.

 \begin{thmx}[\cite{B89,JL24,JL25_1,JL25_2}]\label{leB}
When $d\ge 3$ and $\beta\in [0,\beta_c]$ the limit $W^{\beta}_\infty:= \lim_{n\to \infty} W^{\beta}_n$ is well defined and satisfies $\bbP( W^{\beta}_{\infty}> 0)=1$.
Furthermore for any $p\in [1,1+\frac{2}{d})$, $W^{\beta}_n$ converges to $W^{\beta}_{\infty}$ in $L^p$. Additionally there exists a constant $C>0$ such that
\begin{equation}\label{crititail}
 \forall u\ge 1,\quad   \frac{1}{C} u^{-1-\frac{2}{d}} \le  \bbP\left(  W^{\beta_c}_\infty \ge u \right)
   \le \bbP\left(  \sup_{n\ge 0} W^{\beta_c}_n \ge u \right)\le  C u^{-1-\frac{2}{d}}.
\end{equation}

\end{thmx}
The observation that $W^{\beta}_n$ is a martingale and thus converges as $n\to \infty$ comes from \cite{B89}. It was also proved in the same reference that $\bbP( W^{\beta}_{\infty}> 0)=1$ for small values of $\beta$. The fact that $\bbP( W^{\beta}_{\infty}> 0)=1$ whenever $\f(\beta)=0$ is the main result in  \cite{JL24,JL25_2} (for upper-bounded and general disorder respectivelly). The $L^p$ integrability is a consequence of \cite[Theorem C item (ii)]{J_JEMS}. Finally \eqref{crititail} is a combination \cite[Theorem 2.1]{JL25_1} with \cite[Corollary 2.5]{JL25_2}.

\subsection{Considerations on criticality}\label{conscrit}
When $\beta_c<\infty$ we introduce upper and lower critical exponents, which quantify the regularity of $\f(\beta)$ around $\beta_c$. They are defined as
\begin{equation}\label{nu}
 \overline\nu:= \limsup_{\beta\downarrow \beta_c} \frac{\log |\f(\beta)|}{\log (\beta-\beta_c)} \quad \text{ and } \quad  \underline \nu:=
 \liminf_{\beta\downarrow \beta_c} \frac{\log |\f(\beta)|}{\log (\beta-\beta_c)}.
\end{equation}
They satisfy $2\le \underline \nu \le \overline \nu\le \infty$ (the lower bound is a consequence of the convexity bound on $\f(\beta)+\gl(\beta)$) and we expect that $\underline \nu= \overline \nu$. Their -- hypothetical -- common value, denoted by $\nu$ is referred to as the \textit{critical exponent} associated with $\f$. While the value $\nu$ might depend on the dimension $d$, it is expected that it is independent of the particular choice of distribution for $\go$.

\medskip

The value of $\nu$ has been computed in a couple of cases.
For instance, it has has been shown that $\nu=4$ when $d=1$ and $\nu=\infty$ for $d=2$ (see \cite{L10}, we also refer to \cite{N19} and \cite{BL15} for the sharpest known estimates for $\f(\beta)$ when $\beta\downarrow 0$ in that case) and that $\nu=2$ for the directed polymer on a tree \cite{bpp93}. When $d=1,2$ the value of $\beta_c$ is known to be zero.  On the tree the value of $\f(\beta)$ can be computed for all $\beta>0$, and thus $\beta_c$ is also explicit.

\medskip

When $d\ge 3$, on the contrary,
while some methods have been developped to obtain bounds and approximations for $\beta_c$  (see for instance \cite{CV06,DE92}), there is no established conjecture for the exact value of $\beta_c$. On the contrary the current consensus is that there exists no simple expression for  $\beta_c$ in terms of the distribution of $\go$. An explicit lower bound is given by
$$\beta_2\coloneqq \sup\left\{\beta \ \colon e^{\gl(2\beta)-2\gl(\beta)} E^{\otimes 2}\left[\sum_{n=1}^{\infty} \ind_{\{X^{(1)}_n=X^{(2)}_n\}}\right]< 1\right\}\in(0,\infty].$$
The inequality $\beta_c\ge \beta_2$ can be obtained by the explicit computation of $\bbE[(W^{\beta}_n)^2]$, but it has been proved that this $L^2$ criterion is not sharp. We refer to \cite{Bir04} and \cite[Section 1.4]{BS10} for a proof that $\beta_c>\beta_2$ when $d\ge 4$ and to \cite{BT10} and \cite[Theorem B]{JL24} for the case $d=3$.

\medskip

Not knowing the critical point makes the study of the phase transition utterly challenging, even at the heuristic level. We could find only a couple of references on the topic  in the physics literature. 
Early numerical studies on the topic for $d=3$ yielded conflicting predictions $\nu= 2.1 \pm 0.1$ in \cite{DG90} and $\nu =4\pm 0.7$ in \cite{KBM91}. More recent claims can be found in \cite{MG2,MG1}. In \cite{MG2} it is predicted that $\nu=\infty$ AND $\beta_c=\beta_2$ (and as discussed above, the second part of the claim has since been rigorously disproved) while  \cite{MG1} provides numerical simulations which give two different estimates for the value of $\nu$.\footnote{In the above mentioned references,  \textit{specific heat} and  \textit{correlation lengths} exponents are considered rather than $\nu$ defined by \eqref{nu}, but this difference should not matter as far as heuristics and predictions are concerned.}.

\subsection{Our contribution}

Our main result is a rigorous proof that $\nu=\infty$ for every $d\ge 3$ and every distribution satisfying \eqref{allmoments}.
Our approach relies on studying the growth rate of the moments of the partition function  $\bbE[(W^{\beta}_n)^p]$, in the limit when  $\beta\downarrow \beta_c$, for values of $p$ smaller than  but close to $1+(2/d)$. While idea of using noninteger moments to estimate the free energy dates back to the early years of the study of the model (see for instance \cite{DE92}), the key novelty of our approach is to use conditional expectations and convex comparisons to obtain a bound on $\bbE[(W^{\beta}_n)^p]$ in terms of moments of point-to-point partition functions at $\beta_c$.
This strategy bases itself on the recent understanding of fine integrability  properties of $W^{\beta}_n$ in the subcritical phase and at criticality obtained in a recent series of paper \cite{FJ23,J20,J22,J_JEMS,JL25_1}.

\subsection{Free energy and trajectory behavior}
To conclude this introduction, we provide an illustration of the link between the  free energy and the behavior of the trajectories  under the polymer measure.
It has been shown that the two regimes $\bbP( W^{\beta}_{\infty}> 0)=1$ and $\f(\beta)<0$ correspond to drastically different behavior of the trajectory $(X_k)^n_{k=0}$ under $P^{\beta}_n$. Convergence to Brownian Motion under diffusive rescaling has been proved to hold when $\beta\in [0,\beta_c]$ \cite{CY06} (see also \cite{L25} for a recent concise proof) while endpoint localization has been established when $\beta>\beta_c$ \cite{CSY03} (see also \cite{BC20} for the current most detailed results concerning end-point localization).
There also  exists  quantitative connections between:
\begin{itemize}
 \item [(a)] $\f(\beta)$ and end-point localization,
 \item [(b)] the derivative $\f'(\beta)$ and the asymptotic replica-overlap fraction.
\end{itemize}
To illustrate this, we introduce two random variables which quantify how much trajectories are localized under the polymer measure
\begin{equation}\begin{split}
\mathrm{EP}^{\beta,\go}_n := \frac{1}{n} \sum^n_{k=1}(P^{\beta,\go}_{k})^{\otimes 2}(\ind_{\{X^{(1)}_{k}=X^{(2)}_{k}\} })\\
\mathrm{OV}^{\beta,\go}_n := \frac{1}{n} \sum^n_{k=1}(P^{\beta,\go}_{n})^{\otimes 2}(\ind_{\{X^{(1)}_{k}=X^{(2)}_{k}\} })
\end{split}\end{equation}
The quantity $\mathrm{EP}^{\beta,\go}_n$ is the Cesaro mean of the probability that two polymer replicas share the same end-point while $\mathrm{OV}^{\beta,\go}_n$ corresponds to the overlap fraction between two replicas.

 \begin{thmx}[\cite{CH02,CSY03}]\label{leC}
There exists positive \textbf{continous} functions $\beta\mapsto m(\beta)$ and $\beta \mapsto M(\beta)$ such that for every $\beta>0$ we have almost surely
 \begin{equation}\label{endpointloca}
    -m(\beta) \f(\beta) \le  \liminf_{n\to \infty} \mathrm{EP}^{\beta,\go}_n \le \limsup_{n\to \infty} \mathrm{EP
    }^{\beta,\go}_n \le -M(\beta) \f(\beta).
\end{equation}
If the environment is standard Gaussian ($\go_{1,0}\sim \cN(0,1)$)  then we have
\begin{equation}\label{overlaploca}
 -\f'(\beta+)\le  \beta \liminf_{n\to \infty} \bbE\left[ \mathrm{OV}_n^{\beta,\go} \right]\le \beta \limsup_{n\to \infty} \bbE\left[ \mathrm{OV}^{\beta,\go}_n \right]\le- \f'(\beta-)
\end{equation}
where $\f'(\beta-)$ and $\f'(\beta+)$ denote the left and right derivative of $\f$. Since $\f(\beta)+\gl(\beta)$ is convex  they exist for all $\beta$ and coincide for all except possibly countably many $\beta$s.
\end{thmx}
Since the result does not appear in this exact form in \cite{CH02} and \cite{CSY03}, let us shortly explain how \eqref{endpointloca} and \eqref{overlaploca} can be extracted from the content of these articles.

\medskip

The comparison \eqref{endpointloca} is derived from \cite[Theorem 2.1]{CSY03}, with the additional observation that the constants $c_1$ and $c_2$ appearing in \cite[Equation (2.3)]{CSY03} can be chosen to depend continuously on $\beta$ (this can be checked by going through the proof of \cite[Lemma 3.1]{CSY03}). The continuity in $\beta$ is important when $d\ge 3$ since it implies in that case, $m$ and $M$ can be treated as constant in a neighborhood of $\beta_c$.

\medskip

The comparison in \eqref{overlaploca} is obtained by taking the limit in the identity in \cite[Lemma 7.1]{CH02}. More precisely use the fact that  $(1/n)\log Z^{\beta,\go}_n$ is convex in $\beta$ to interchange limit and derivative (see e.g.\ \cite[Section A.1.1]{GiacBook} for comments on this  interchange of limits).

\section{Results}

\subsection{Smoothness of the free energy and consequences}
Our main result states that $f(\beta)$ vanishes faster than any power of $(\beta-\beta_c)$ at criticality.
\begin{theorem}\label{supermain}
For the directed polymer on $\bbZ^d$ with $d\ge 3$ we have
 \begin{equation}\label{supersmooth}
   \lim_{\beta \downarrow \beta_c }\frac{\log |\f(\beta)|}{\log (\beta-\beta_c)}=\infty.
 \end{equation}
\end{theorem}

As a consequence of Theorem \ref{leC}, we obtain assymptotic estimates for  $\mathrm{EP}^{\beta,\go}_n$ and $\bbE\left[ \mathrm{OV}^{\beta,\go}_n \right]$.

\begin{corollary}\label{EPOV}
For every $K>0$ we have
\begin{equation}\label{epsmooth}
 \lim_{\beta  \downarrow \beta_c} (\beta-\beta_c)^{-K} \limsup_{n\to \infty} \mathrm{EP
    }^{\beta,\go}_n =0.
\end{equation}
If the environment is Gaussian ($\go_{1,0}\sim \cN(0,1)$)  then for every $K>0$,
\begin{equation}\label{oversmooth}
  \lim_{\beta  \downarrow \beta_c}  (\beta-\beta_c)^{-K}\limsup_{n\to \infty} \bbE\left[ \mathrm{OV}^{\beta,\go}_n \right]=0.
\end{equation}
\end{corollary}


\subsection{Critical asymptotic behavior for moments of order $p$}

The smoothness result presented above is also valid for some of the Lyapunov exponents associated
with fractional moments of $W^{\beta}_n$.
Given $p> 1$, we set
\begin{equation}\label{lyapudef}
\tf(p,\beta):= \lim_{n\to \infty}   \frac{1}{n}\log \bbE\left[  \left(W^{\beta}_n\right)^p\right]=\lim_{n\to \infty}\frac 1 n \log \bbE[(Z^{\beta,\go}_n)^p]-p\gl(\beta).
\end{equation}
As shown in  \cite[Proposition 2.1]{J22} the limit exists by superadditivity since for every $m\in \bbN$ and $n\in \bbN\cup\{\infty\}$ we have
\begin{equation}\label{ddddo}
\bbE\left[  \left(W^{\beta}_{m+n}\right)^p \ | \ \cF_m \right]\le   \left(W^{\beta}_{m}\right)^p \bbE\left[ \left(W^{\beta}_n\right)^p\right]
\end{equation}
and hence $\bbE\left[  \left(W^{\beta}_{m+n}\right)^p \right] \le \bbE\left[ \left(W^{\beta}_{m}\right)^p\right] \bbE\left[ \left(W^{\beta}_n\right)^p\right]$.
Since we are going to make use of \eqref{ddddo} in our proof (more precisely in  Section \ref{qdtr}) we include a proof  at the end of this introduction for completeness. Let us display here some basic properties of  $\tf(p,\beta)$ as a function of $\beta$.

\begin{lemma}\label{bazprop}
The function   $\beta\mapsto \tf(p,\beta)$ is continuous nondecreasing.
Hence given $p>1$ there exists $\beta_p$ such that
 \begin{equation}\label{lesbetac}
  \tf(p,\beta)>0    \quad \Leftrightarrow \quad \beta> \beta_p \ .
 \end{equation}
When $p\in (1,1+\frac{2}{d}]$, we have $\beta_p=\beta_c$.
\end{lemma}

\begin{proof}[Sketch of proof]
The monotonicity in $\beta$ of $\tf(p,\beta)$ is obtained by using convex comparison (see Section \ref{convsec}) and is a direct consequence of Lemma \ref{monoconvex}. One way to prove continuity is to observe that $\tf(p,\beta)+\gl(\beta)$ is convex in $\beta$ by checking via a computation (recall \eqref{lyapudef} above) that  $\partial^2_{\beta}\log \bbE[(Z^{\beta,\go}_n)^p]\ge 0$.
The inequality $\beta_p\ge  \beta_c$ for  $p\in (1,1+\frac{2}{d}]$ is a consequence of  Theorem \ref{leB}, and $\beta_p\le \beta_c$ for every $p>1$ can be deduced from Lemma  \ref{comparisionprim} presented below.
\end{proof}

We show that, when $p\in (1,1+\frac{2}{d})$ the phase transition for the $L^p$ moments is also smooth.

\begin{theorem}\label{themomentboundprim}
For any $p\in (1,1+\frac{2}{d})$ we have
  \begin{equation}
   \lim_{u\downarrow 0}\frac{\log \tf(p,\beta_c+u)}{\log u}=\infty.
 \end{equation}
\end{theorem}

\subsection{Near critical behavior}

As a by product of the techniques used in the proof of Theorem \ref{themomentboundprim}, we obtain information on the \textit{near critical} behavior of the directed polymer. If $\beta=\beta_n$ is allowed to depend on $n$ and $\beta_n-\beta_c$ converges to $0$ like an inverse power of $n$ with an \textit{arbitrarily small} exponent then the limiting partition function is $W^{\beta_c}_\infty$ \footnote{The use of $\beta_n$ here is in slight conflict with that of $\beta_p$ in Lemma \ref{bazprop}, but this should not create any serious confusion.}.

\begin{proposition}\label{sideresult}

Given $\delta>0$ and a sequence $(\beta_n)_{n\ge 1}$ such that 
\begin{equation}
 \forall n\ge 1, \beta_n>\beta_c \quad \text{  and  } \quad  \lim_{n \to \infty} n^{\delta} (\beta_n-\beta_c)= 0,
\end{equation}
then we have $\lim_{n\to \infty} W^{\beta_n}_{n}= W^{\beta_c}_\infty$
where the convergence holds in $L^p$ for all $p\in (1,1+\frac{2}{d})$.
\end{proposition}

\subsection{Open questions and further research directions}

\subsubsection*{Free energy asymptotics around $\beta_c$}
A natural question that comes to mind when reading Theorem \ref{supermain} is ``Can one find a simple equivalent of $\f$ in the neighborhood of $\beta_c$?''.
In analogy with the critical behavior for the directed polymer for  $d=2$ \cite{BL15} and similar estimates obtained for a continuum polymer model with a  Poisson environment \cite[Theorem 2.11, item (ii)]{BCL23} (valid for any $d\ge 2$)  we may expect that there exists $\gamma(d)\in (0,\infty)$ such that
\begin{equation}\label{conjecturexx}
 \f(\beta)\stackrel{\beta\downarrow \beta_c}{\sim} \exp\left(- (\beta-\beta_c)^{-\gamma(d)+o(1)}\right).
\end{equation}
In the current state of affairs obtaining \text{any} upper bound on $\f(\beta)$ that goes beyond \eqref{supersmooth} or any lower bound, presents a significant technical challenge.

\subsubsection*{Intermediate regime}
To understand further the transition from weak to strong disorder, and one may look for an \textit{intermediate disorder}  regime in which   $\beta_n\downarrow \beta_c$ as $n\to \infty$ and  $W^{\beta_n}_{n}$ converges (in some sense) to a limit which is nontrivial and different
from $W^{\beta_c}_\infty$ (see \cite{AKQ14,CSZ17,CSZ23}  and references for corresponding results in dimension $1$ and $2$). In view of Proposition \ref{sideresult}, we know that to obtain such a regime, $\beta_n-\beta_c$ should be taken larger than any negative power of $n$, but identifying the exact scaling is a completely open question.

\subsubsection*{Regularity of the free energy}
Since $\beta_c$ marks the only known phase transition for the directed polymer model, one may expect that $\f(\beta)$ is analytic on $\bbR \setminus \{\beta_c\}$ and
 \eqref{supersmooth} further indicates that  $\f(\beta)$ should be infinitely differentiable at $\beta_c$.
However our proof method does not provide any tool to control the regularity of $\f$ at $\beta_c$ or elsewhere.

%
%

\subsection{Organization of the paper}

In Section \ref{conexion}, we present a comparison result between $\f(\beta)$ and $\tf(p,\beta)$. As a  consequence of this comparison, Theorem \ref{supermain} can be deduced from Theorem \ref{themomentboundprim}. The comparison is obtained by exploiting the concentration of $\log W^{\beta}_n$ around its mean, established in  \cite{LW09, W10}.

\medskip

In Section  \ref{proofprim}, we expose the main idea behind the proof of Theorem  \ref{themomentboundprim}. We also provide a decomposition of this proof which relies of two main technical results: Propositions \ref{compapinningz} and \ref{comparisionprim}.
These propositions are proved in Section \ref{lalaya} and \ref{qdtr} respectively.

\medskip

Section \ref{oxysec} is dedicated to the proof of auxiliary results:  Corollary \ref{EPOV} and Proposition \ref{sideresult}.

\subsection{Notation}
We close this section by introducing all the notation we use in the remainder of the paper.
 We let $\cF_n:= \sigma\left( \go_{k,x} :   k\in \lint 1, n \rint\  , \ x\in \bbZ^d  \right)$  (we use the notation $\lint a,b\rint:= [a,b]\cap \bbZ$) define the natural filtration $(\cF_n)_{n\ge 0}$ associated with $\go$. We recall (this can be checked by hand) that $(W^{\beta}_n)_{n\ge 1}$ is a martingale w.r.t.\ with this filtration.
 In particular $\bbE\left[ (W^{\beta}_n)^p\right]$ is increasing in $n$ when $p\ge 1$.
We define the \textit{point-to-point} partition function by setting
\begin{equation}
\hat W^{\beta}_n(x)=E \left[ e^{\sum_{j=1}^n(\beta \go_{j,X_j}-\gl(\beta))}\ind_{\{X_n=x\}}\right].
\end{equation}
For $k\ge 0$ and $z\in \bbZ^d$ the shift operator $\theta_{k,z}$ on $\go$ is defined by
\begin{equation}\label{definishift}
 \theta_{k,z} \go:= (\go_{n+k,x+z})_{n\ge 1,x\in \bbZ^d}.
\end{equation}
We let $\theta_{k,z}$ act on functions of $\go$ by setting $\theta_{k,z} f=f\circ\theta_{k,z}$ (we use this mainly for partition functions). Finally we introduce the multiplicative weights $\zeta_\beta$ setting
\begin{equation}\label{zetadef}
\zeta_{\beta}(k,x):=e^{\beta\go_{k,x}-\gl(\beta)} \quad \text { and } \quad  \zeta_{\beta}:= \zeta_{\beta}(1,0).
\end{equation}
We assume for convenience and without loss of generality that $\bbE[\go_{1,0}]=0$ and $\bbE[\go^2_{1,0}]>0$ this implies that
\begin{equation}\label{stand}
\forall \beta>0, \ \gl(\beta)>0 \quad \text{ and } \quad  \gl'(\beta)>0.
\end{equation}
The notation being introduced we can present a proof a statement presented earlier in this introduction.
\begin{proof}[Proof of \eqref{ddddo}]
With the notation introduced above, we observe that by decomposing the trajectory according to the value of $X_m$ we have
\begin{equation}
W^{\beta}_{m+n}=\sum_{x\in \bbZ^d} \hat W^{\beta}_m(x)\theta_{m,x}W^{\beta}_{n}.
\end{equation}
Hence using Jensen's inequality for the probability $\hat W_m(\cdot)/W_m$ we obtain
\begin{equation}
 (W^{\beta}_{m+n})^p\le (W^{\beta}_m)^{p} \sum_{x\in \bbZ^d} \frac{\hat W^{\beta}_m(x)}{W_m}(\theta_{m,x}W^{\beta}_{n})^p
\end{equation}
and finally using translation invariance and independence we obtain
\begin{equation}\begin{split}
  \bbE\left[ (W^{\beta}_{m+n})^p \ | \ \cF_m \right]&\le  (W^{\beta}_m)^{p} \sum_{x\in \bbZ^d} \frac{\hat W^{\beta}_m(x)}{W_m}\bbE\left[ (\theta_{m,x}W^{\beta}_{n})^p \ | \ \cF_m \right]\\
  & = (W^{\beta}_m)^{p}  \bbE\left[ (\theta_{m,x}W^{\beta}_{n})^p\right]. \qedhere
  \end{split}
\end{equation}
\end{proof}

\section{Connection between the free energy and the moment growth}\label{conexion}

\subsection{The main statement}
\label{compargu}

The main achievement of this section is a comparison between the free energy and and $\tf(p,\cdot)$ around the critical point. The comparison is valid for all $p>1$ but only relevant for $p\in (1, 1+\frac{2}{d}]$ which are the values of $p$ for which $\beta_p=\beta_c$ (see Lemma \ref{bazprop}).

\begin{proposition}\label{comparisionprim}
We set $K_1:= 2 e^{\gl(\beta_c+1)+ \gl(-\beta_c-1)}$.
Given $p>1$ for all $u\in [0,1]$ such that $|\f(\beta_c+u)|\le 2K_1$ we have
 \begin{equation}\label{pasmaal}
 |\f(\beta_c+u)|\le  2 \sqrt{ \frac{K_1}{p-1} \tf(p,\beta_c+u)}.
 \end{equation}
 Hence for any $p\in (1,1+\frac{2}{d}]$, we have
\begin{equation}\label{exposantin}
    \liminf_{u\downarrow 0}\frac{\log |\f(\beta_c+u)|}{\log u}\ge \frac{1}{2} \liminf_{u\to 0+}    \frac{\log \tf(p,\beta_c+u)}{\log u}.
\end{equation}
\end{proposition}

  Note that in view of \eqref{exposantin}, Theorem \ref{supermain}  follows from  Theorem \ref{themomentboundprim}. In the remainder of the section, we present two technical results required for the proof, before proving Proposition \ref{comparisionprim}.

\subsection{Concentration inequality}

A consequence of the assumption  \eqref{allmoments} is that the increments of the  Doob martingale $\big(M^{(n)}_{k}\big)_{k=0}^m:=\big(\bbE\left[ \log W^{\beta}_n \ | \ \cF_k \right]\big)^n_{k=0}$ have uniformly bounded \textit{conditional} exponential  moments $\bbE\left[  e^{\gl(M_{k}-M_{k-1})} \ | \ \cF_{k-1}\right]$. This property implies exponential concentration around the mean. We refer to  \cite{LW09,W10} for detailed results. We are going to make use of the following estimate.

\begin{proposition}\cite[Theorem 1.2]{W10} \label{fromW}
For $K= 2 e^{\gl(\beta)+ \gl(-\beta)}$ and any $x\ge 0$ we have
\begin{equation}
      \bbP\left( \log W^{\beta}_n- \bbE\left[\log W^{\beta}_n \right]\ge x\right) \le e^{-n \varphi_K(x/n)} \quad \text{ where }
      \quad \varphi_K(u):=\begin{cases}
                            \frac{u^2}{4K}, \text{ if } u\le K, \\
                            u-K \text{ if } u\ge K.
                           \end{cases}
\end{equation}

\end{proposition}

\subsection{A variant of Payley-Zygmund inequality}
This is the $L^p$ variant of the classical Paley-Zygmund inequality (together with its proof).
\begin{lemma}\label{PZ}
Given $Z$ a nonnegative random variable with mean one and $p>1$ we have
\begin{equation}
 \bbP\left[Z \ge \theta\right] \ge  (1-\theta)^\frac{p}{p-1}\bbE\left[ Z^p\right]^{-\frac{1}{p-1}}.
\end{equation}
\end{lemma}
\begin{proof}
 This is a direct consequence of H\"older's inequality applied as follows
 \begin{equation}
 1=\bbE[Z]\ge \bbE[Z\ind_{\{Z\ge \theta\}}]+\theta \ge  \bbP(Z\ge \theta)^{\frac{p-1}{p}} \bbE\left[ Z^p\right]^{\frac{1}{p}}+\theta. \qedhere
 \end{equation}
\end{proof}

\subsection{Proof of Proposition \ref{comparisionprim}}
Using Lemma \ref{PZ} we have
\begin{equation}
 \bbP\left[ W^{\beta_c+u}_n \ge 1/2\right] \ge  2^{-\frac{p}{p-1}}\bbE\left[ (W^{\beta_c+u}_{n})^{p}\right]^{-\frac{1}{p-1}}.
\end{equation}
and hence
\begin{equation}\label{le10}
 \liminf_{n\to \infty}\frac{1}{n}\log  \bbP\left( W^{\beta_c+u}_n \ge 1/2\right)\ge - \frac{1}{p-1} \tf(p,\beta_c+u).
\end{equation}
On the other hand for $u\in [0,1]$,
we have from Proposition \ref{fromW} (and the fact that $\varphi_K$ is decreasing in $K$)
\begin{equation}
  \bbP\left( W^{\beta_c+u}_n \ge 1/2\right) \le  \exp\left(- n\varphi_{K_1}\left(-\frac{\log 2}{n}  -\frac{1}{n}\bbE \left[ \log  W_n\right]  \right) \right)
\end{equation}
Taking $\log$, dividing by $n$ and letting $n$ tend to infinity one obtains
\begin{equation}\label{le20}
  \limsup_{n\to \infty}\frac{1}{n}\log  \bbP\left( W^{\beta_c+u}_n \ge 1/2\right)\le -\varphi_{K_1}\left( \f(\beta_c+u) \right)= -  \frac{\f(\beta_c+u)^2}{4K_1}
\end{equation}
where the last equality is valid provided that  $|\f(\beta_c+u)|< 2K_1$.
The inequality \eqref{pasmaal} follows from the combination \eqref{le10} and \eqref{le20}.
Then \eqref{exposantin} follows since by continuity the condition $\f(\beta_c+u)\le 2K_1$ is satisfied for $u$ sufficiently small.
\qed

\section{Proof of Theorem  \ref{themomentboundprim}}\label{proofprim}

\subsection{The starting idea: conditional second moment}\label{startee}
Let us open this section by exposing one of the main underlying idea for the proof. 
At the starting point, we look for an upper bound for $|\f(\beta)|$ rather than for $\tf(p,\beta)$.

\medskip

  Both for $d=1$ and $d=2$, reasonably precise upper bounds on $|\f(\beta)|$ have been obtained (in the $\beta\to 0$ limit) by estimating (or computing) the growth rate of  $\bbE\left[ (W^{\beta}_n)^2 \right]$ for small values of $\beta$,  and combine the estimate with either a percolation argument (see for instance \cite[Section 5]{L10} or the proof of  \cite[Theorem 1.3]{AY15}) or concentration estimates (e.g.\ as in Section \ref{conexion} above and in \cite[Section 7]{L10}). This is not possible to adopt such a strategy for $d\ge 3$ because, since $\beta_c>\beta_2$ (recall the content of Section \ref{conscrit}), the second moment sequence  $\bbE\left[ (W^{\beta}_n)^2 \right]$ starts to blow up exponentially before $\beta$ reaches criticality.

\medskip

What we do instead is to consider a \textit{conditional} second moment.
Let us focus on the Gaussian case ($\go_{1,0}\sim \cN(0,1)$) for which the idea is simpler to expose
 (Proposition \ref{convexdom} displays  how the reasoning can be adapted to general environment).
The starting observation is that for $\beta>\beta_c$, the \textit{supercritical} Gibbs weights $\zeta_{\beta}(k,x)$ (recall \eqref{zetadef}) have the same distribution as  the product of \textit{critical} Gibbs weights $\zeta_{\beta_c}(k,x)$  with  independent random variables whose variance are of order $\beta-\beta_c$. This is  due to the following equality in distribution
\begin{equation}\label{specicoupling}\begin{split}
\left((\beta_c+u) \    \go_{k,x} \right)_{k\ge 1, x\in \bbZ^d} &\stackrel{(d)}{=} \left(\beta_c \   \go_{k,x}+ \  \sqrt{v} \ \go'_{k,x} \right)_{n\ge 1,x\in \bbZ^d},\\
 \left(\zeta_{\beta_c+u}(k,x)\right)_{k\ge 1, x\in \bbZ^d} &\stackrel{(d)}{=} \left(\zeta_{\beta_c}(k,x) e^{\sqrt{v}\go'_{k,x}-\frac{v}{2}}\right)_{k\ge 1, x\in \bbZ^d}
\end{split}\end{equation}
where $\go$ and $\go'$ are two indendent i.i.d.\ standard Gaussian fields (with distribution $\bbP$ and $\bbP'$) and $v= 2\beta_c u+ u^2$.
We set
\begin{equation}
 \overline W^{\beta_c+u}_n:= E\left[ e^{\sum_{k=1}^n \left(\beta_c \go_{k,X_k}+ \sqrt{v} \go'_{k,X_k}- \frac{\beta^2_c+v}{2} \right)}\right]
 = W^{\beta_c}_n\times  E^{\beta_c,\go}_n\left[  e^{\sum_{k=1}^n \left(\sqrt{v} \go'_{k,X_k}- \frac{v}{2}\right) }\right]
\end{equation}
From \eqref{specicoupling}, the distribution of  $\overline W^{\beta_c+u}_n$ under $\bbP\otimes \bbP'$ is equal to that of  $W^{\beta_c+u}_n$.
We use Fubini to compute the second moment of  $\overline W^{\beta_c+u}_n$ with respect to $\go'$ and we obtain
\begin{equation}\begin{split}\label{toop}
 \bbE'\left[  \left(\overline W^{\beta_c+u}_n\right)^2\right]
 &=  (W^{\beta_c}_n)^2 (E^{\beta_c,\go}_n)^{\otimes 2}\left[ \bbE'\left[ \exp\left(\sum_{k=1}^n \left[ \sqrt{v} (\go'_{k,X^{(1)}_k}+\go'_{k,X^{(2)}_k}) - v\right] \right) \right]\right]\\
 &= (W^{\beta_c}_n)^2 (E^{\beta_c,\go}_n)^{\otimes 2}\left[ e^{ \sum_{k=1}^n v \ind_{\{X^{(1)}_k=X^{(2)}_k\}}}\right].
\end{split}\end{equation}
Using a conditional analogue of Proposition \ref{comparisionprim} (and $W^{\beta_c}_\infty>0$) we should have 
\begin{equation}\label{condprop31}
 |\f(\beta_c+u)|\le C \sqrt{\limsup_{n\to \infty} \frac{1}{n} \log  (E^{\beta_c,\go}_n)^{\otimes 2}\left[ e^{ \sum_{k=1}^n v \ind_{\{X^{(1)}_k=X^{(2)}_k\}}}\right] }
\end{equation}
and the control of the free energy reduces to that of the growth of the Laplace transform of the replica overlap at criticality.

\medskip

The best solution we found to control the growth rate of this conditional second moment is to bound its moment of order $p/2$
\begin{equation}\label{orderp2}
\bbE\left[ \bbE'\left[(\overline W^{\beta_c+u}_n)^{2}\right]^{p/2} \right]= \bbE\left[ (W^{\beta_c}_n)^{p} (E^{\beta_c,\go}_n)^{\otimes 2}\left[ e^{ \sum_{k=1}^n v \ind_{\{X^{(1)}_n=X^{(2)}_n\}}}\right]^{p/2} \right]
\end{equation}
for values of $p$ smaller but close to $1+\frac{2}{d}$. This approach allows to exploit the main information we have about the directed polymer at $\beta_c$, which the tail the exponent $(1+\frac{2}{d})$ for the asymptotic distribution of the partition function (see Theorem \ref{leB}).
By Jensen's inequality  $\bbE\left[ \bbE'\left[(\overline W^{\beta_c+u}_n)^{2}\right]^{p/2} \right]\le \bbE\otimes \bbE'\left[(\overline W^{\beta_c+u}_n)^{p} \right]$ and thus this reasoning also yields a  proof of
Theorem \ref{themomentboundprim}.

\subsection{Two key estimates}
We now present the technical details of the proof of Theorem \ref{themomentboundprim}, which is obtained by combining two comparison results.
 The first is an upper bound on  $\bbE\left[ \left(W^{\beta_c+u}_n\right)^{p}\right]$ which is obtained -- in the Gaussian case -- by combining an expansion of the conditionnal second moment displayed in \eqref{toop} with the subadditivity of  $s\mapsto s^{p/2}$. To display it we introduce the quantity
\begin{equation}
\mathcal K^{(p)}(n):= \bbE\left[ \sum_{y\in \bbZ^d} \hat W^{\beta}_{n} (y)^{p} \right].
\end{equation}

\begin{proposition}\label{compapinningz}
Given $p\in[1,2]$, there exist constants $C,u_0>0$, which  may both  depend on $\beta$ and on the distribution of $\go$, such that for all $u\in[0,u_0]$, ($i_0=0$ by convention)
\begin{equation}\label{cbien}\begin{split}
\frac{ \bbE\left[ \left(W^{\beta_c+u}_n\right)^{p}\right]}{\bbE\left[ (W^{\beta_c}_{n})^{p}\right]} &\le \cZ_n(Cu^{p/2})\\
\text{ where  } \quad  \cZ_n(w):=
1+ \sum_{k\ge 1} &\sum_{0<i_1<\dots<i_k\le n} w^{k} \prod_{j=1}^k\cK^{(p)}(i_{j}-i_{j-1}).\end{split}
\end{equation}
 
\end{proposition}

The proof of Propoposition \ref{compapinningz} is provided in Section \ref{lalaya}.
The quantity $\cZ_n(w)$ is analogous to the partition function of a free-boundary homogenous pinning model (see \cite[Chapter 2]{GiacBook} for a review on that subject), the main difference being that $\cK^{(p)}(\cdot)$ is not summable. The exponential growth rate of $ \cZ_n(w)$ can thus be controlled by replicating the proof of \cite[Theorem 2.2]{GiacBook}. We show in the next section that  for any $w>0$  (the inequality below is an equality but we will not use this fact  in our proof)
\begin{equation}\label{dropp}
 \limsup_{n\to \infty} \frac{1}{n}  \log \cZ_n(w)\le \varphi(w) \quad \text{ where }  \quad \sum_{n\ge 1} e^{-n\varphi(w)}\cK^{(p)}(n)= \frac{1}{w}.
\end{equation}
To estimate $\varphi(w)$ we need a control over the following quantity
\begin{equation}\label{defjnp}
 J_N(p):= \sum_{n=1}^N \cK^{(p)}(n)= \sum^N_{n=1} \sum_{x\in \lint -N,N\rint^d} \bbE\left[ \left(\hat W^{\beta_c}_n(x)\right)^{p}\right].
 \end{equation}
This is our second key estimate.

\begin{proposition}\label{crucialprim}
 We have, for any $p\in (1,1+\frac{2}{d})$
  \begin{equation}\label{p}
\limsup_{N\to \infty} \frac{ \log J_N(p) }{\log N}\le d+2-dp.
 \end{equation}
 in particular we have 
 \begin{equation}\label{onlyneed}
\lim_{p\uparrow (1+\frac{2}{d})}\limsup_{N\to \infty}  \frac{\log J_N(p)}{\log N} =0.
\end{equation}
\end{proposition}
The proof of Proposition \ref{crucialprim} is presented in Section \ref{qdtr}

\begin{remark}
 Note that by Jensen's inequality for $p \in (1,1+\frac{2}{d})$ we have
 \begin{equation}\label{yuio}
   J_N(p)\ge \sum^N_{n=1} \sum_{x\in \lint-N,N \rint^d} P(X_n=x)^{p} = (C_{d,p}+o(1)) N^{\frac{d+2-dp}{2}}.
 \end{equation}
Adding some technical refinements in the proof of Proposition \ref{crucialprim} one may prove that this lower bound is sharp as far as the exponent is concerned in the sense that
   \begin{equation}\label{pprim}
\lim_{N\to \infty} \frac{ \log J_N(p) }{\log N}= \frac{d+2-dp}{2}.
 \end{equation}
Since for Theorem \ref{themomentboundprim}, we only need \eqref{onlyneed},
we only prove the easier \textit{ad-hoc} bound \eqref{p}.
\end{remark}

\subsection{The actual proof of Theorem \ref{themomentboundprim}}

First we observe that it is sufficient to show that 
 for every $p\in (1,1+\frac{2}{d})$
\begin{equation}\label{alternative}
\liminf_{u\downarrow 0} \frac{\log \tf(p,\beta_c+u)}{\log u}\ge \frac{p}{2(d+2-dp)}.
\end{equation}
Indeed by Jensen's inequality
$p\mapsto  \frac{1}{p}\tf(p,\beta)$ is nondecreasing and hence if $1\le p\le q< 1+\frac{2}{d}$, \eqref{alternative} implies that
\begin{equation}\label{alternative2}
\liminf_{u\downarrow 0} \frac{\log \tf(p,\beta_c+u)}{\log u} \ge \liminf_{u\downarrow 0} \frac{\log \tf(q,\beta_c+u)}{\log u} \ge \frac{q}{2(d+2-dq)},
\end{equation}
and we conclude by letting $q$ tend to $1+\frac{2}{d}$.

\medskip

Let us now prove \eqref{alternative}. We start with a proof of \eqref{dropp}.  Setting $\tilde K(n):= w e^{-n\varphi(w)} \cK^{(p)}(n)$, we let $\tau$ (with law $\bP$) be a random integer sequence with $\tau_0=0$ and whose increments $(\tau_i-\tau_{i-1})_{i\ge 1}$ are i.i.d.\ with distribution given $\tilde K$  we have
\begin{equation}\begin{split}\label{daup}
\cZ_{n}(w)&= 1+ \sum_{k \ge 1} \sum_{0<i_1<\dots<i_k\le  n} e^{\varphi(w) i_k} \prod_{j=1}^k\tilde \cK(i_{j}-i_{j-1})\\
&\le e^{\varphi(w)n} \left( 1+ \sum_{k \ge 1} \sum_{0<i_1<\dots<i_k\le n}  \prod_{j=1}^k\tilde \cK(i_{j}-i_{j-1}) \right)\\
&= e^{\varphi(w)n}\left(  \sum_{m=0}^n \bP\left(  \exists k\ge 0,  \tau_k=m  \right) \right)\le (n+1)e^{\varphi(w)n}
\end{split}\end{equation}
As a consequence combining \eqref{cbien} with \eqref{dropp}, we have for $p\in (1,1+\frac{2}{d})$ (for these values of $p$, $\bbE\left[ (W^{\beta_c}_n)^p \right]$ is bounded by Theorem \ref{leB})
\begin{equation}
 \tf(p,\beta_c+u)=\lim_{n\to \infty} \frac{1}{n}\log \bbE\left[ (W^{\beta_c+u}_n)^p \right]
 \le   \limsup_{n\to \infty} \frac{1}{n}\log \cZ_{n}(Cu^{p/2}) \le \varphi(C u^{p/2}).
\end{equation}
As a consequence we have
\begin{equation}
\liminf_{u\downarrow 0} \frac{\log \tf(p,\beta_c+u)}{\log u}\ge \frac{p}{2}\liminf_{w\downarrow 0} \frac{\log \varphi(w)}{\log w}
\end{equation}
We conclude our proof of \eqref{alternative} by showing that
\begin{equation}\label{lait}
\liminf_{w\downarrow 0} \frac{\log \varphi(w)}{\log w}\ge \frac{1}{d+2-dp}.
\end{equation}

Since  $J_N(p):=\sum_{n=1}^N \cK^{(p)}(n)$, using summation by part, we obtain that
\begin{equation}
  \frac{1}{w}  = \sum_{n\ge 1} e^{-n\varphi(w)}\cK^{(p)}(n)=(1-e^{-\varphi(w)}) \sum_{N\ge 1}e^{-\varphi(w)N} J_N(p).
\end{equation}
From Proposition \ref{crucialprim}, given $\delta>0$ there exits positive constants $C_{\delta}$ and $C'_{\delta}$ such that for all $v\in [0,\varphi^{-1}(1)]$ we have
\begin{equation}
 (1-e^{-\varphi(w)}) \sum_{N\ge 1} J_N(p)e^{-N\varphi(w) }\le   C_{\delta}\varphi(w)\sum_{N\ge 1} N^{d+2-dp+\delta} e^{-\varphi(w) N} \le C'_{\delta} \varphi(w)^{dp-d+2-\delta}.
\end{equation}
We obtain that $\varphi(w) \le \left(C'_{\delta} w\right)^{\frac 1 {d+2-dp+\delta}},$ which concludes the proof of \eqref{lait} since $\delta$ is arbitrary.
\qed

\section{Proof of Proposition \ref{compapinningz}}\label{lalaya}

\subsection{Proof of Proposition \ref{compapinningz} in the Gaussian case}

 Recalling the content and notation of Section \ref{startee} (and in particular  that  $v=2\beta_c u+ u^2$ is of the same order as $u$) we obtain from  \eqref{orderp2} that
 \begin{equation}\label{gotit}
   \bbE\left[ (W^{\beta_c+u}_n)^p\right] = \bbE\otimes \bbE'\left[ (\overline W^{\beta_c+u}_n)^p\right]
   \le  \bbE\left[ (W^{\beta_c}_n)^{p} (E^{\beta_c,\go}_n)^{\otimes 2}\left[ e^{ \sum_{k=1}^n v \ind_{\{X^{(1)}_k=X^{(2)}_k\}}}\right]^{p/2} \right]
\end{equation}

Then, setting $\delta_I:=\ind_{\{\forall k\in I, \  X^{(1)}_k=X^{(2)}_k\}}$ we expand the exponential as follows
$$\exp\Big(u\sum_{k=1}^n \ind_{\{X^{(1)}_k=X^{(2)}_k\}}\Big)=\prod_{k=1}^n (1+ (e^v-1)\ind_{\{X^{(1)}_k=X^{(2)}_k\}})= \sum_{I\subset \lint 1,n\rint}  \left(e^v-1\right)^{|I|}  \delta_I.$$
Pushing the decomposition further, we let  $(i_k)^{|I|}_{k=1}$ denote the ordered enumeration of the elements of $I$ and for ${\bf x}=(x_1,\dots, x_{|I|})\in (\bbZ^{d})^{|I|}$ we set
$\delta(I,{\bf x}):=\ind_{\{\forall k\in \lint 1, |I|\rint, \  X^{(1)}_{i_k}=X^{(2)}_{i_k}=x_{k}\}}.$ Taking expection w.r.t.\ $(P^{\beta_c,\go}_n)^{\otimes 2}$ we obtain
\begin{equation}
 (E^{\beta_c,\go}_n)^{\otimes 2}\left[ e^{ \sum_{k=1}^n v \ind_{\{X^{(1)}_k=X^{(2)}_k\}}}\right]= \sum_{I\subset \lint 1,n\rint}  \sum_{{\bf x}\in (\bbZ^d)^{|I|}} \left(e^v-1\right)^{|I|}
 (E^{\beta_c,\go}_n)^{\otimes 2} \left[\delta (I,{\bf x})\right].
\end{equation}
Combining this with subadditivity (recall that by assumption $p/2\le 1$), multiplying by $(W^{\beta_c}_n)^{p}$ and taking expectation w.r.t.\ $\bbP$, we obtain finally (recall \eqref{gotit})
\begin{equation}\begin{split}\label{gotitt}
\bbE\left[ (W^{\beta_c+u}_n)^p\right]&\le  \bbE\left[  (W^{\beta_c}_n)^{p} (E^{\beta_c,\go}_n)^{\otimes 2}\left[ e^{ \sum_{k=1}^n v \ind_{\{X^{(1)}_k=X^{(2)}_k\}}}\right]^{p/2}\right]\\
  &\le  \sum_{I\subset \lint 1,n\rint}  \sum_{{\bf x}\in (\bbZ^d)^{|I|}} (e^v-1)^{\frac{2|I|}{p}}
    \bbE\left[ \left( (W^{\beta_c}_n)^2(E^{\beta_c,\go}_n)^{\otimes 2} \left[\delta (I,{\bf x})\right]\right)^{p/2} \right].
\end{split}\end{equation}
We observe that the expression $(W^{\beta_c}_n)^2(E^{\beta_c,\go}_n)^{\otimes 2}[\delta (I,{\bf x})]$ is a product of squared partition functions. Using the convention $i_0=0$ and $x_0=0$,
 we have (recall \eqref{definishift})
\begin{equation}
(W^{\beta_c}_n)^2(E^{\beta_c,\go}_n)^{\otimes 2}[\delta (I,{\bf x})]
= \left(\prod_{k=1}^{|I|} \theta_{i_{k-1},x_{k-1}} \hat W^{\beta_c}_{i_k-i_{k-1}}(x_{k}-x_{k-1})\right)^2 \left(\theta_{i_{|I|},x_{|I|}} W^{\beta_c}_{n-x_{|I|}}\right)^2.
\end{equation}
Using independence and translation invariance after taking the above equality to the power $p/2$, we obtain that
\begin{equation}
\bbE\left[ \left( (W^{\beta_c}_n)^2(E^{\beta_c,\go}_n)^{\otimes 2} \left[\delta (I,{\bf x})\right]\right)^{p/2} \right] =\prod_{k=1}^{|I|} \bbE\left[ \left(\hat W^{\beta_c}_{i_k-i_{k-1}}(x_{k}-x_{k-1})\right)^{p}\right] \bbE\left[ \left(W^{\beta_c}_{n-x_{|I|}}\right)^{p}\right].
\end{equation}
 We use the above in the r.h.s.\ of \eqref{gotitt}. We re-index the sum, setting $(y_k)^{|I|}_{k=1}:=(x_{k}-x_{k-1})_{k=1}^{|I|}$ and we obtain
\begin{equation}\begin{split}
   \bbE\left[ \left(W^{\beta_c+u}_N\right)^{p} \right]
   &\le  \sum_{I\subset \lint 1,n\rint} \left(e^v-1\right)^{\frac
{|I| p}{ 2}  }  \bbE\left[ \left(W^{\beta_c}_{n-i_{|I|}}\right)^{p}\right]\prod_{k=1}^{|I|} \sum_{y_k\in \bbZ^d }\bbE\left[ \left(\hat W^{\beta_c}_{i_k-i_{k-1}}(y_{k})\right)^{p}\right] \\
&\le \sum_{I\subset \lint 1,n\rint} \left(e^v-1\right)^{\frac
{|I| p}{ 2}  }  \bbE\left[ \left(W^{\beta_c}_{n}\right)^{p}\right]\prod_{k=1}^{|I|} \cK^{(p)}(i_{k}-i_{k-1}) ,
\end{split}\end{equation}
where  the last inequality simply comes from the fact that   $\bbE\left[ \left(W^{\beta_c}_{n}\right)^{p}\right]$ is increasing in $n$.
Since $e^v-1\le C u$ for all $u\le 1$ (recall $v=2\beta_c u+u^2$) we obtain the desired result.
\qed

\subsection{Convex comparison}\label{convsec}
The proof of Proposition \ref{compapinningz} in the general case relies on the notion of convex comparison to which we now provide a short introduction (see \cite[Chapter 3]{SS07} for a more thorough review).
Given two random variables in $Z_1$ and $Z_2$ in $L^1$, we say that $Z_2$ convexly dominates $Z_1$ and write
$Z_1 {\preccurlyeq}_\mathrm{(conv)} Z_2$ if  $\bbE[Z_1]=\bbE[Z_2]$ and for all convex function $\varphi: \bbR\to \bbR_+$ we have
\begin{equation}\label{defconvorder}
 \bbE\left[ \varphi(Z_1) \right] \le \bbE\left[ \varphi(Z_2)\right].
\end{equation}
The terminology is slightly improper since ${\preccurlyeq}_\mathrm{(conv)}$ defines a partial order on probability distributions rather than on $L^1$. Approximating $\varphi$ by functions which are affine by part, one can check (see \cite[Theorem 3.A.1]{SS07}) that if $\bbE[Z_1]=\bbE[Z_2]$ then \eqref{defconvorder} is equivalent to the following ($x_+:=\max(x,0)$)
\begin{equation}\label{onlypospart}
\forall a \in \bbR, \quad  \bbE\left[ (Z_1-a)_+ \right]\le  \bbE\left[ (Z_2-a)_+ \right].
\end{equation}
With an induction  on $M$  we obtain the following extension of  \eqref{defconvorder}.
\begin{lemma}\label{muilvariate}
Let $M\ge 1$,  $(Y_1,Z_1),\dots,(Y_M,Z_M)$ be independent $\bbR^2$-valued random variables with $Y_i \preccurlyeq_\mathrm{(conv)} Z_i$ for every $i\in \lint 1, M\rint$,
and $\phi: \bbR^M\to \bbR_+$ be such that for every $i\in \lint 1, M\rint$ and $x\in \bbR^M$,`
$  u\mapsto \phi(x_1,\dots,x_{i-1},u, x_{i+1},\dots, x_M)$ is convex,
 then $$\bbE\left[ \phi(Y_1,\dots,Y_M)\right]\le\bbE\left[ \phi(Z_1,\dots,Z_M)\right].$$
\end{lemma}
The following result illustrates how the notion of convex order can be useful in the context of directed polymers (recall \eqref{zetadef}).
\begin{lemma}\label{monoconvex}
 The sequence $(\zeta_{\beta})_{\beta\ge 0}$ is monotone increasing for the convex order.
 As a consequence, the sequence
 $(W^{\beta}_n)_{\beta\ge 0}$ is monotone increasing for the convex order.
\end{lemma}

The result was proved first as \cite[Lemma 3.3]{CY06}, we include a short proof for completeness.

\begin{proof}
 By a limit argument, it is only necessary to prove that $\bbE\left[\varphi(\zeta_{\beta})\right]$ is monotone in $\beta$ for convex functions $\varphi$ that are both Lisphitz and $C^1$. Recalling that from the definition in \eqref{zetadef}, we have $\partial_{\beta} \zeta_\beta= (\go_{1,0}- \gl'(\beta)) \zeta_{\beta}$, we obtain
 \begin{equation}
  \partial_{\beta}\bbE\left[ \varphi(\zeta_{\beta}) \right]=  \bbE\left[ \partial_{\beta}\varphi(\zeta_{\beta}) \right]= \bbE\left[ (\beta \go_{1,0}-\gl'(\beta))  \zeta_{\beta}\varphi'(\zeta_{\beta}) \right].
 \end{equation}
 Since $\zeta_\beta$ is positive and $\bbE[\zeta_{\beta}]=1$, we may define a probability $\tilde \bbP$ by setting $\dd \tilde \bbP= \zeta_\beta \dd \bbP$.
 Note that both $(\beta \go_{1,0}-\gl'(\beta))$ and $\varphi'(\zeta_{\beta})$ are increasing functions of $\go_{1,0}$ ($\varphi$ is convex).
 Since  $\int f g \  \dd \mu \ge \int f \dd \mu \int g \dd \mu$ for any pair of increasing $L^2(\mu)$ functions and any probability $\mu$
we have
$$\tilde \bbE\left[ (\beta \go_{1,0}-\gl'(\beta)) \varphi'(\zeta_{\beta}) \right]\ge \tilde \bbE\left[  (\beta \go_{1,0}-\gl'(\beta)) \right]\tilde \bbE\left[ \varphi'(\zeta_{\beta})\right].$$
Finally, we  note that  $\tilde \bbE\left[  (\beta \go_{1,0}-\gl'(\beta)) \right]=\bbE\left[ (\beta \go_{1,0}-\gl'(\beta)) \zeta_{\beta}\right]=0$ to conclude. To establish the result concerning the partition function we apply Lemma \ref{muilvariate} after observing that $\varphi(W^{\beta}_{n})$ is a convex function of $\zeta_{\beta}(k,x)$ for every $(k,x)\in \lint 1,n\rint\times \lint -n,n\rint^d$.
\end{proof}

\subsection{Proof of Proposition \ref{compapinningz} in the general case}\label{lalaya2}

In the above proof, the Gaussian environment assumption is used in \eqref{specicoupling} to provide a useful martingale coupling of $\zeta_{\beta_c}(k,x)$ and $\zeta_{\beta_c+u}(k,x)$. While the identity \eqref{specicoupling} is very specific to Gaussian environment,
convex comparisons provide an alternative possibility.

\begin{proposition}\label{convexdom}
There exists $C_{\beta}>0$ and $u_0(\beta)>0$ such that for all $u\in [0, u_0]$ we can construct a nonnegative random variable $Y_{u}=Y^{(\beta)}_u$ independent of $\zeta_{\beta}$  which satisfies the following
 \begin{equation}
\bbE[Y_u]=1, \quad   \zeta_{\beta+u}{\preccurlyeq}_{\mathrm{(conv)}} Y_u \zeta_{\beta} \quad \text{ and } \quad
\bbE\left[  Y_u^2 \right]  \le \left( 1+       C_{\beta} u\right).
 \end{equation}
\end{proposition}
The proof of Proposition \ref{convexdom} is postponed to Section \ref{domdom}.
We consider $(Y^{(u)}_{k,x})_{k\ge 1, x\in \bbZ^d}$ to be i.i.d. random variables, independent of $\go$ satisfying
\begin{itemize}
 \item [(i)] $\bbE\left[Y^{(u)}_{k,x}\right]=1$ \text{ and }  $\bbE[(Y^{(u)}_{k,x})^2]\le 1+Cu$,
  \item [(ii)]   $e^{(\beta_c+u)\go_{nkx}-\gl(\beta_c+u)}{\preccurlyeq}_{\mathrm{(conv)}}Y^{(u)}_{n,x} e^{\beta_c\go_{k,x}-\gl(\beta_c)}  $.
\end{itemize}
We let $\bbP'$ denote the distribution of  the field $Y^{(u)}$. Using
Lemma \ref{muilvariate} we obtain that
\begin{equation}
W^{\beta_c+u}_n {\preccurlyeq}_{\mathrm{(conv)}}
E\left[ \prod_{k=1}^n  \left(e^{\beta_c \go_{k,X_k}- \gl(\beta_c)}Y^{(u)}_{k,X_k}\right)\right]
=W^{\beta_c}_n E^{\beta_c,\go}_{n}\left[ \prod_{k=1}^n Y^{(u)}_{k,X_k}\right].
\end{equation}
Using the  convex comparison (first inequality), Jensen's inequality (the second inequality) and finally Fubini  like in \eqref{toop} together with item $(i)$ (third inequality) we obtain
\begin{equation}\begin{split}\label{debutprim}
 \bbE\left[ \left(W^{\beta_c+u}_n\right)^{p} \right]
  &\le \bbE\otimes \bbE'\left[ (W^{\beta_c}_n)^p E^{\beta_c,\go}_{n}\left[ \prod_{k=1}^n Y^{(u)}_{k,X_k}\right]^{p} \right]\\
  &\le \bbE\left[   (W^{\beta_c}_n)^p \bbE' \left[ E^{\beta_c,\go}_{n}\left[ \prod_{k=1}^n Y^{(u)}_{k,X_k}\right]^{2} \right]^{\frac p 2} \right]\\
  &\le \bbE\left[   (W^{\beta_c}_n)^p (E^{\beta_c,\go}_{n})^{\otimes 2}\left[ (1+Cu\ind_{\{ X^{(1)}_{k}=X^{(2)}_k\}}) \right]^{p/2}\right].
 \end{split}
\end{equation}
Then we can then repeat the computation performed for the Gaussian case with  $(e^v-1)$ replaced by $Cu$. \qed

 \subsection{Proof of Proposition \ref{convexdom}}\label{domdom}
For $v\in[0,1/3]$ we let $Z_v$ be a random variable which is independent of $\zeta_{\beta}$ and distributed as follows
\begin{equation}\begin{cases}
 \bbP(Z_v\in \dd  r)=  6v r^{-4}\dd r \quad \text{ for } r\in (1,\infty) , \\ \bbP(Z_{v}=1)=1- 3 v,  \\  \bbP(Z_v=0)=v.
\end{cases}\end{equation}
By construction $\bbE[Z_v]=1$ and $\bbE[Z_v^2]=1+3v.$
We are going to show that there exists a constant $C=C_{\beta}$ such that for all $u\in[0, 1]$ we have $\zeta_{\beta+u} {\preccurlyeq}_{\mathrm{(conv)}}  \zeta_{\beta}Z_{Cu}$.
Recalling \eqref{onlypospart}, it is sufficient to show that  following holds for $u$ sufficiently small
\begin{equation}\label{comparesidz}
\forall x\in \bbR, \quad
 \bbE\left[ (\zeta_{\beta+u}-x)_+-(\zeta_{\beta}-x)_+\right]\le \bbE\left[ (\zeta_{\beta}Z_{Cu} -x)_+-(\zeta_{\beta}-x)_+\right].
\end{equation}
Since both sides are equal to zero if $x\le 0$ we only need to adress the case $x>0$.
Using the notation $a_+=\max(a,0)$ and letting $\partial^+_{\beta}$ denote the derivative on the right we have
 \begin{equation}\label{ripp}\begin{split}
  \partial^+_{\beta}\bbE\left[ (\zeta_{\beta}-x)_+\right]&= \bbE\left[ (\go-\gl'(\beta))\zeta_{\beta}\ind_{\{\zeta_{\beta}>x\}}+  (\go-\gl'(\beta))_+  \zeta_{\beta}\ind_{\{\zeta_{\beta}=x\}} \right]\\
  &= \bbE\left[ (\gl'(\beta)-\go)\zeta_{\beta}\ind_{\{\zeta_{\beta}<x\}}+ (\gl'(\beta)-\go)_+\zeta_{\beta}\ind_{\{\zeta_{\beta}=x\}} \right]
 \end{split}\end{equation}
 (the second line is obtained simply by observing that $\bbE\left[ (\go-\gl'(\beta))\zeta_{\beta}\right]=0$).
Integrating the first line \eqref{ripp} on the interval $[\beta,\beta+u]$, we obtain for any $u\in[0, 1]$
 \begin{equation}\begin{split}\label{small3}
 \bbE\left[ (\zeta_{\beta+u}-x)_+- (\zeta_{\beta}-x)_+\right] &\le \int^u_0\bbE\left[ (\go-\gl'(\beta+v))_+ \zeta_{\beta+v}\ind_{\{\zeta_{\beta+v}\ge x\}}\right] \dd v  \\ &\le u \bbE\left[ \go_+ e^{(\beta+1)\go_+}\ind_{\{e^{(\beta+1)\go_+}\ge x\}}\right]\\
 & \le   u \bbE\left[\go^2_+ e^{2(\beta+1)\go_+} \right]^{1/2} \bbP\left(\go_+ \ge \frac{\log x}{\beta+1} \right)^{1/2}.
\end{split}\end{equation}
In the second line we used \eqref{stand} and the last inequality is simply Cauchy-Schwarz.
We need another bound for $x$ small, using the second line in \eqref{ripp} we have
 \begin{equation}\begin{split}\label{small2}
 \bbE\left[ (\zeta_{\beta+u}-x)_+- (\zeta_{\beta}-x)_+\right] &\le\int^u_0\bbE\left[ (\gl'(\beta+v)-\go)_+\zeta_{\beta+v}\ind_{\{\zeta_{\beta+v}\le x\}}\right] \dd v  \\ &\le u x \left(\gl'(\beta+1)+ \bbE\left[|\go|\right]\right).
\end{split}\end{equation}
The combination of \eqref{small3} and \eqref{small2} yields that there exists $C$ such that that all $u\in [0,1]$ and all $x\in(0,\infty)$ we have
\begin{equation}\label{small4}
 u^{-1} \bbE\left[ (\zeta_{\beta+u}-x)_+- (\zeta_{\beta}-x)_+\right]\le C \min\left(x, \bbP\left(\go_+ \ge \frac{\log x}{\beta+1} \right)^{1/2}\right).
\end{equation}

On the other hand the r.h.s.\ in \eqref{comparesidz} has a simple explicit expresssion.
\begin{equation}\begin{split}\label{small1}
 v^{-1}\bbE\left[ (\zeta_{\beta}Z_v-x)_+-   (\zeta_{\beta}-x)_+\right]
 &= 6\int^{\infty}_1 r^{-4}  \bbE\left[ (r\zeta_{\beta}-x)_+\right] \dd r   -  3\bbE\left[ (\zeta_{\beta}-x)_+\right]\\
 &= x\bbP(\zeta_{\beta}\ge x)+ 6\int^{\infty}_1 r^{-4}  \bbE\left[ (r\zeta_{\beta}-x)_+\ind_{\{\zeta_\beta<x\}}\right]\dd r.
\end{split}\end{equation}
The second line in \eqref{small1} is obtained by observing that
\begin{multline}
6\int^{\infty}_1 r^{-4}  \bbE\left[ (r\zeta_{\beta}-x)_+\ind_{\{\zeta_{\beta}\ge x\}} \right] \dd r=  6\int^{\infty}_1 \left(r^{-3} \bbE\left[ \zeta_{\beta}\ind_{\zeta_{\beta}\ge x\}} \right]+ r^{-4} x\bbP(\zeta_{\beta}\ge x)  \right)\dd r      \\ =   3\bbE[\zeta_{\beta}\ind_{\{\zeta_{\beta}\ge x\}} ]
-2 x \bbP(\zeta_{\beta}\ge x)
=3\bbE\left[ (\zeta_{\beta}-x)_+\right]  + x \bbP[\zeta_{\beta}\ge x].
\end{multline}
In view of \eqref{small4} and \eqref{small1}, to conclude our proof  it is sufficient to show that $\sup_{x>0} g(x) < \infty$ where
\begin{equation}\label{finalcompa}
 g(x):=\frac{\min\left(x, \bbP(\go_+ \ge \frac{\log x}{\beta+1} )^{1/2}\right)}{ \left(x\bbP(\zeta_{\beta}\ge x) +   \int^{\infty}_1 r^{-4}  \bbE\left[ (r\zeta_{\beta}-x)_+\ind_{\{\zeta_\beta<x\}}\right]\dd r \right)}.
\end{equation}
Since $g$ is countinuous, it is sufficient to check that it is bounded in neighborhoods of $0$ and infinity. Since $g(x)\le 1/\bbP(\zeta_{\beta}\ge x)$
we have
\begin{equation}
\limsup_{x\downarrow 0} g(x)\le 1.
\end{equation}
For large values of $x$ we observe that
\begin{equation}
  \bbE\left[ (r\zeta_{\beta}-x)_+\ind_{\{\zeta_\beta<x\}}\right]\ge x\bbP\left( \zeta_{\beta}\in(1/2,x) \right)\ind_{\{r\ge 4x\}},
\end{equation}
which integrated over $r\in (1,\infty)$ yields
\begin{equation}
 \int^{\infty}_1 r^{-4}  \bbE\left[ (r\zeta_{\beta}-x)_+\ind_{\{\zeta_\beta<x\}}\right]\dd r\ge  \frac{\bbP\left( \zeta_\beta\in( 1/2,x)\right)}{192 x^2}.
\end{equation}
Since \eqref{allmoments} implies that $\bbP(\go_+ \ge \frac{\log x}{\beta+1} )$ decays faster than any power of $x$, we have
\begin{equation}
 \limsup_{x\to \infty} g(x)\le   \limsup_{x\to \infty} \frac{\bbP(\go_+ \ge \frac{\log x}{\beta+1} )^{1/2}}{ \int^{\infty}_1 r^{-4}  \bbE\left[ (r\zeta_{\beta}-x)_+\ind_{\{\zeta_\beta<x\}}\right]\dd r}=0,
\end{equation}
concluding the proof of \eqref{finalcompa}.
\qed

\section{Proof of Proposition \ref{crucialprim}}\label{qdtr}

\subsection{The starting idea}
 Let us briefly expose the strategy of proof. We know from Theorem \ref{leB} that $\sup_{M\ge 0}\bbE[(W^{\beta_c}_{\infty})^p]<\infty$ is finite.
 We want to ``compare'' $J_N(p)$ to $\bbE[(W^{\beta_c}_{\infty})^p]$.
 By restricting the partition function to trajectories satisfying  $X_n=x$ we observe that
 \begin{equation}\label{simple}
 \forall n \ge 1, \quad  \forall x\in \bbZ^d    \quad  \quad  W^{\beta_c}_\infty \ge  \hat W^{\beta_c}_n(x)  \theta_{n,x} W^{\beta_c}_{\infty}.
 \end{equation}
We want to exploit the above bound in subtle way.
We consider a family of events $(\cA_{n,x})_{n\in \lint 1, N\rint,x\in \lint -N,N\rint^d}$ which satisfy
\begin{itemize} 
 \item [(a)] The event are disjoint:  $\cA_{n,x}\cap \cA_{n',x'}=\emptyset$ if $(n,x)\ne (n',x')$,
 \item [(b)] $\cA_{n,x}$ is independent of $\hat W^{\beta_c}_n(x)$ for every $x$ and $n$.
\end{itemize}
From \eqref{simple} and disjointness we have 
\begin{equation}
 (W^{\beta_c}_\infty)^p\ge \sum_{n\in \lint 1, N\rint}\sum_{x\in \lint -N,N\rint^d} (\hat W^{\beta_c}_n(x))^p  \left(\theta_{n,x} W^{\beta_c}_{\infty}\right)^p \ind_{\cA_{n,x}}.
\end{equation}
Taking expectation and using the independence property we have 
\begin{equation}\begin{split}\label{simplez}
 \bbE\left[ (W^{\beta_c}_\infty)^p  \right]& \ge
 \sum_{n\in \lint 1, N\rint}\sum_{x\in \lint -N,N\rint^d} \bbE\left[(\hat W^{\beta_c}_n(x))^p\right]  \bbE\left[\left(\theta_{n,x} W^{\beta_c}_{\infty}\right)^p \ind_{\cA_{n,x}}\right]\\
 &\ge  J_{N}(p) \min_{(n,x)\in \lint 1, N \rint \times \lint -N,N\rint^d}\bbE\left[\left(\theta_{n,x} W^{\beta_c}_{\infty}\right)^p \ind_{\cA_{n,x}}\right].
\end{split}\end{equation}
The challenge is  to find a partition $\cA_{n,x}$ which is such that the minimum above is not too small. A natural idea is to start with $A^0_{n,x}:= \{ \theta_{n,x} W^{\beta_c}_{\infty}\ge N^K \}$ for some $K>0$ (which satisfies $(b)$ above but not $(a)$) and then turn it into a disjoint family by considering some $\cA_{n,x}\subset A^{0}_{n,x}$ which only cuts off a ``small'' portion of the original event.
From Theorem \ref{leB} there exists constant $c,C>0$ such that
\begin{equation}\label{iouiou}
\bbE\left[\left(W^{\beta_c}_{\infty}\right)^p \ind_{\{W^{\beta_c}_{\infty}\ge N^{K}\}}\right]\ge c N^{K(1+\frac{2}{d}-p)} \quad \text{ and } \quad
\bbP\left( W^{\beta_c}_{\infty}\ge N^{K} \right) \le C N^{-K\left(1+\frac{2}{d} \right)}.
\end{equation}
The first inequality in \eqref{iouiou} indicates that (for any fixed $K$ and $\gep$) we have
$$\bbE\left[\left(\theta_{n,x} W^{\beta_c}_{\infty}\right)^p \ind_{A^{0}_{n,x}}\right]\ge c N^{-\gep}$$ for $p$ sufficiently close to $1+\frac{2}{d}$. This is  what we require for the proof \eqref{onlyneed}. The second inequality implies that for $K> \frac{d(d+1)}{(d+2)}$ (in the next subsection we choose $K=d$)
\begin{equation}\label{whatheK}
\sum_{n\in \lint 1, N\rint,x\in \lint -N,N\rint^d}\bbP(A^{0}_{n,x})\le C'N^{1+d-K\left(1+\frac{2}{d}\right)} \ll 1
\end{equation}
which leaves some hope that $(A^{0}_{n,x})$ can be turned into a disjoint collection of events without chopping a big portion off the events .

\subsection{Its implementation}
To implement the above idea, we need to perfom a couple of modification. Firstly, to in order have a finite range dependence, instead of considering $\{  \theta_{n,x} W^{\beta_c}_{\infty}\ge N^d \}$ we set 
\begin{equation}\label{defanx}
A_{n,x}:= \{ \theta_{n,x}W^{\beta_c}_{m_N}> N^d \}.
\end{equation}
We want to take $m_N$ as small as possible, but in a way that keeps $\bbP\left( W^{\beta_c}_{m_N}> N^d \right)$ roughly comparable to $\bbP\left( W^{\beta_c}_{\infty}> N^d \right)$.
To choose the value of $m_N$ we rely on a result from  \cite{JL25_1} which compares the tail distribution of $W^{\beta_c}_{m}$ and that of  $W^{\beta_c}_{\infty}$ up to  to a certain threshold.
\begin{lemma}[{\cite[Lemma 5.2]{JL25_1}}]\label{neuefunfzwei}
For any $\gep>0$ there exists $C_{\gep}$ and $u_0(\gep)$ such that for all $u\ge u_0$
\begin{equation}\label{ptoline}
 \exists m\in \lint 1,  C_{\gep}\log u\rint, \quad  \bbP\left(  W^{\beta_c}_{m}\ge  u \right)\ge u^{-1-\frac{2}{d}-\gep}.
\end{equation}
\end{lemma}
Let us clarify here that there is no mystery on how to chose $m$ in \eqref{ptoline}.
We set $q= 1+\frac{1}{d}$. Using first Lemma \ref{PZ}  and then  \eqref{ddddo}, we have 
\begin{equation}\begin{split}
\bbP\left(  W^{\beta_c}_{n+m}\ge \frac{1}{2} W^{\beta_c}_{m} \ | \   W^{\beta_c}_{m}\ge 2u   \right) &\ge 2^{-(d+1)} \bbE\left[ (W^{\beta_c}_{n+m}/ W^{\beta_c}_{m})^q \ | \  W^{\beta_c}_{m}\ge u \right]^{-d}\\
\\&\ge  2^{-(d+1)} \bbE\left[ (W^{\beta_c}_{n})^q\right]^{-d}.
\end{split}\end{equation}
This implies that (recall that $W^{\beta_c}_n$ converges in $L^q$, see Theorem \ref{leB})
\begin{equation}
\bbP\left(  W^{\beta_c}_{n+m}\ge u\right)\ge 2^{-(d+1)} \bbE\left[ (W^{\beta_c}_{\infty})^q\right]^{-d}\bbP\left( W^{\beta_c}_{m}\ge 2u  \right).
\end{equation}
Hence 
for a different choice of $C_{\gep}$ and $u_0$, we have
\begin{equation}\label{corolemma}
  \forall m \ge C_{\gep}\log u, \quad  \bbP\left(  W^{\beta_c}_m\ge  u \right)\ge u^{-1-\frac{2}{d}-\gep}.
\end{equation}
We set $m_N:= (\log N)^2$ and fix $\delta>0$.
 Using \eqref{corolemma} with $u=N^{d}$ (for the lower bound) and Theorem \ref{leB} (for the upper bound) we have for $N$ sufficiently large
\begin{equation}\label{fraaming}
 2 N^{-2-d-\delta}\le \bbP\left(  W^{\beta_c}_{m_N} \ge N^d \right) \le C N^{-2-d}.
\end{equation}
where $C$ is the constant given by Theorem \ref{leB}.

\medskip

Secondly, we modify slightly the disjointness assumption. With the definition \eqref{defanx}, since $\ind_{A_{n,x}}$ is a function of $(\go_{k,z})_{k\in \lint n+1, n+m_N\rint,\ z\in x+ \lint -m_N,m_N\rint^d}$,  the collection of event $(A_{n,x})$ has a finite range of dependence. More precisely, we consider the  partition
$$\lint 1,N\rint \times \lint -N,N\rint^d:= \bigsqcup^{M_N}_{j=1} \mathcal J_j$$ where $(\cJ_j)_{j=1}^{M_N}$ is an enumeration of the intersection of $\lint 1,N\rint \times \lint -N,N\rint^d$ with the cosets of $(m_N\bbZ)\times ((2m_N+1)\bbZ)^d$ (with $M_N:= m_N (2m_N+1)^d$). For every $j$ the events $(A_{n,x})_{(n,x)\in \cJ_j}$ are independent.
Since independent collections of event are easier to work with, we are going to treat each coset separetely. For this reason we replace the assumptions considered
of the previous subsection  by ($(b)$ is unchanged)
\begin{itemize} 
 \item [(a$\star$)]  For any $j\in \lint 1, M_N\rint$ the events $(\cA_{n,x})_{(n,x)\in \cJ_j}$  are disjoint, in particular we have 
 $$\sum_{(n,x)\in \lint 1,N\rint \times \lint -N,N\rint^d} \ind_{\cA_{n,x}}\le M_N.$$ 
  \item [(b)] $\cA_{n,x}$ is independent of $\hat W^{\beta_c}_n(x)$ for every $x$ and $n$.
 \end{itemize}
 We are now ready to define $\cA_{n,x}$.
 We equip $\bbN\times \bbZ^d$ with the lexicographical order (we have $(n',x')> (n,x)$ if either $n'>n$ or if $n=n'$ and $x'>x$ for the  lexicographical order on $\bbZ^d$) and for  $(n,x)\in \mathcal J_j$ we set
\begin{equation}\begin{split}
A'_{n,x} &:= A_{n,x}\cap \bigcap_{x'\in \bbZ^d :  (n,x')\in \mathcal J_j \text{ and } x'>x } A^{\complement}_{n,x'}\\
\mathcal B_{n} &:= \bigcap_{(n',x')\in \mathcal J_j :  \text{ and } n'>n } A^{\complement}_{n,x'}\\
\mathcal A_{n,x}&:= A_{n,x} \cap \bigcap_{(n',x')\in \mathcal J_j \ : \ (n',x')> (n,x)} A^{\complement}_{n',x'}= A'_{n,x}\cap \mathcal B_{n}.
\end{split}\end{equation}
In words, $\mathcal A_{n,x}$ is the event ``$(n,x)$ is the largest element of $\mathcal J_j$ such that $\theta_{n,x}W^{\beta_c}_m> N^d$''.
By construction the events $\cA_{n,x}$ satisfy $(a\star)$ and $(b)$. For $(b)$ it is sufficient to observe that 
$\cA_{n,x}\in \sigma( \go_{k,z}, k\ge n+1, z\in \bbZ^d)$.
From  $(a\star)$ we obtain (using the bound \eqref{simple})
\begin{equation}
  (W^{\beta_c}_\infty)^p\ge \frac{1}{M_N}\sum_{n\in \lint 1, N\rint}\sum_{x\in \lint -N,N\rint^d} (\hat W^{\beta_c}_n(x))^p  \left(\theta_{n,x} W^{\beta_c}_{\infty}\right)^p \ind_{\cA_{n,x}}.
\end{equation}
Then using (b) as in \eqref{simplez} we obtain
\begin{equation}\label{topzz}
 J_{N}(p)\le M_N \bbE\left[  (W^{\beta_c}_\infty)^p \right] \left( \min_{(n,x)\in \lint 1, N \rint \times \lint -N,N\rint^d}\bbE\left[\left(\theta_{n,x} W^{\beta_c}_{\infty}\right)^p \ind_{\cA_{n,x}}\right]\right)^{-1}.
\end{equation}
Since $\bbE\left[  (W^{\beta_c}_\infty)^p \right]$ is simply a constant, and $M_N$ grows like a power of $\log N$, to conclude the proof of Proposition \ref{crucialprim} it is sufficient to prove the following.
\begin{lemma}\label{pote}
 For $N\ge N_0(\delta)$ sufficiently large, we have
 \begin{equation}\label{topzzz}
\forall n\in \lint 1, N \rint,\ \forall x\in  \lint -N,N\rint^d, \quad   \bbE\left[\left(\theta_{n,x} W^{\beta_c}_{\infty}\right)^p \ind_{\cA_{n,x}}\right]\ge  \frac 1 2 N^{dp-d-2-\delta}.
 \end{equation}
\end{lemma}
Indeed the the combination of \eqref{topzz} and \eqref{topzzz} yields
\begin{equation}
  \limsup_{N\to \infty} \frac{J_{N}(p)}{\log N}\le dp-d-2-\delta,
\end{equation}
and we conclude by sending $\delta\to 0$. \qed

\subsection{Proof of Lemma \ref{pote}}

From \eqref{fraaming} we obtain that for any $N$ sufficiently large and $(n,x)\in \mathcal J_j$
($|\cJ_j|$ is of order $(N/m_N)^{d+1}\ll N^{d+1}$)
\begin{equation}\begin{split}\label{ofinterest}
 \bbP(\cB^{\complement}_n)&\le  \sum_{(n',x')\in  \mathcal J_j}\bbP(A_{n',x'})\le C|\mathcal J_j| N^{-d-2}\le N^{-1},\\
  \bbP(A'_{n,x})&= \bbP(A_{n,x}) \prod_{x'\in \bbZ^d \  : \  (n,x')\in \mathcal J_j \text{ and } x'>x }\bbP\left( A^{\complement}_{n,x'}\right)\ge N^{-d-2-\delta}
 \end{split}\end{equation}
where the inequalities are a consequence of \eqref{fraaming}. For the second line, we use the upper bound in \eqref{fraaming} to show that $\prod_{x'\in \bbZ^d : \cdots} \cdots \ge 1/2$ for $N$ large enough.
Since $\cA_{n,x}:=  A'_{n,x}\cap \cB_n$ and $A'_{n,x}$ is $\cF_{n+m}$ measurable we have
\begin{equation}\begin{split}\label{you1}
 \bbE\left[  (\theta_{n,x}W^{\beta_c}_{\infty})^p \ind_{\cA_{n,x}} \ | \ \cF_{n+m} \right]&= (\theta_{n,x} W^{\beta_c}_m)^p \ind_{A'_{n,x}}\bbE\left[ \theta_{n,x} \left(\frac{W^{\beta_c}_{\infty}}{W^{\beta_c}_m}\right)^p\!\!\! \ind_{\cB_{n}} \ | \ \cF_{n+m}\right]\\
 \end{split}\end{equation}
 To estimate the last term we rewrite it as
 \begin{equation}
 \bbE\left[ \theta_{n,x} \left(\frac{W^{\beta_c}_{\infty}}{W^{\beta_c}_m}\right)^p\ | \ \cF_{n+m}\right]- \bbE\left[ \theta_{n,x} \left(\frac{W^{\beta_c}_{\infty}}{W^{\beta_c}_m}\right)^p\!\!\! \ind_{\cB_{n}^{\complement}} \ | \ \cF_{n+m}\right]
 \end{equation}
Using Jensen's inequality for the conditional expectation, we have
\begin{equation}\label{you2}
 \bbE\left[ \theta_{n,x} \left(\frac{W^{\beta_c}_{\infty}}{W^{\beta_c}_m}\right)^p\ | \ \cF_{n+m}\right]\ge 1.
\end{equation}
On the other hand, setting $p'= \frac{d+2+2p}{2d}$ (the midpoint between $p$ and $1+\frac{2}{d})$) and using first  H\"older's inequality for the conditional expectation, and then \eqref{ddddo} (combined with translation invariance) and \eqref{ofinterest} we obtain that for $N$ sufficiently large
\begin{equation}\begin{split}\label{you3}
\bbE\left[ \theta_{n,x} \left(\frac{W^{\beta_c}_{\infty}}{W^{\beta_c}_m}\right)^p\!\!\! \ind_{\cB_{n}^{\complement}} \ | \ \cF_{n+m}\right]&\le  \bbE\left[ \theta_{n,x} \left(\frac{W^{\beta_c}_{\infty}}{W^{\beta_c}_m}\right)^{p'} \ | \ \cF_{n+m}\right]^{\frac{p}{p'}} \bbP(\cB^{\complement}_n)^{\frac{p'-p}{p'}}\\
&\le  \bbE\left[ (W^{\beta_c}_\infty)^{p'} \right]^{\frac{p}{p'}}  N^{-\frac{p'-p}{p'}}\le  1/2.
\end{split}
\end{equation}
The finiteness of $\bbE\left[ (W^{\beta_c}_\infty)^{p'} \right]$ is ensured by Theorem \ref{leB}. Combining \eqref{you1}, \eqref{you2} and \eqref{you3} and the fact that
$\theta_{m,x} W^{\beta_c}_m\ge N^d$ on $A'_{n,x}$ we obtain
\begin{equation}
  \bbE\left[  (\theta_{n,x}W^{\beta_c}_{\infty})^p \ind_{\cA_{n,x}} \ | \ \cF_{n+m} \right]\ge  \frac{1}{2} N^{dp} \ind_{A'_{n,x}}.
\end{equation}
Taking the expectation and using \eqref{ofinterest}, we obtain that
\begin{equation}  \bbE\left[  (\theta_{n,x}W^{\beta_c}_{\infty})^p \ind_{\cA_{n,x}} \right] \ge  \frac{1}{2} N^{dp} \bbP(A'_{n,x})\ge \frac{1}{2}N^{-dp-d-2-\delta}.
\end{equation}
\qed

 \section{Proof of auxiliary results}\label{oxysec}
 
\subsection{Proof of Corollary \ref{EPOV}}
\label{secEPOV}

Setting $\overline{M}:=\max_{\beta\in [\beta_c,\beta_c+1]} M(\beta)$, then for every $\beta\in [\beta_c,\beta_c+1]$, Theorem \ref{leC} implies that
\begin{equation}
\limsup_{n\to \infty} \mathrm{EP}^{\beta,\go}_n\le \overline M \f(\beta)
\end{equation}
holds almost surely, and thus \eqref{epsmooth} is a direct consequence of \eqref{supersmooth}.

\medskip

For \eqref{oversmooth}, we use the convexity of $\f(\beta)+\frac{\beta^2}{2}$ and we  obtain that for any $u>0$
\begin{equation}
 \f'(\beta-)+\beta\ge \frac{\beta^2-(\beta-u)^2+  2(\f(\beta)-\f(\beta-u))}{2u}.
\end{equation}
This yields
\begin{equation}
  -\f'(\beta-) \le   \frac{u}{2}+ \frac{\f(\beta-u)-\f(\beta)}{u}\le  \frac{u}{2}- \frac{\f(\beta)}{u}
  \le \frac{3}{2} \sqrt{|\f(\beta)|}.
\end{equation}
where the last inequality is obtained by taking $u= \sqrt{|\f(\beta)|}$. Hence \eqref{oversmooth}  follows by combining Theorem \ref{leC} and Theorem \ref{supermain}.

\subsection{A martingale estimate}

The following \textit{ad hoc} estimate allows to control the $p$-moment of a martingale difference and is important in the proof of Proposition \ref{sideresult}.

\begin{lemma}\label{otptilem}
 Let $p\in(1,2)$ and $X$ and $Y$ be two $L^p$-integrable random variables such that $\bbE[Y | X]=0$ then if $\bbE\left[ |Y+X|^p-|X|^p \right]\le \bbE[|X|^p]$ we have

 \begin{equation}\label{dfty}
 \bbE[|Y|^p]\le C_p \bbE[|X|^p]^{1-\frac{p}{2}}\bbE[ |Y+X|^p-|X|^p]^{p/2}
 \end{equation}
where $C_p>0$ is a universal constant that depends only on $p$.

\end{lemma}

\begin{proof}
First we observe that for  $p\in(1,2)$, there exists $c_p>0$ such that for all $a,b\in \bbR$
\begin{equation}\label{dsuy}
 |a+b|^p-|a|^p-p\frac{|a|^p}{a}b \ge c_p\min\left(|a|^{p-2} b^2, |b|^{p}\right).
\end{equation}
By homogeneity, it is sufficient to verify that  $|1+b|^p- 1- p b \ge c_p\min\left( b^2, |b|^{p}\right).$\\ 
That can be done by checking separately that 
\begin{equation}
 \begin{cases}
  |1+b|^p- 1- p b\le c^{(p)}_1 b^2 \quad  &\text{ for } |b| \le 1,\\
  |1+b|^p- 1- p b\le  c^{(p)}_2 b^p \quad &\text{ for } |b| \ge 1.
 \end{cases}
 \end{equation}
 (we may take $c^{(p)}_1=p$ and $c^{(p)}_2=2^p$ but this is of no importance for the application we have in mind).
From our assumptions we deduce that
$ \frac{|X|^p}{X}Y$ is in $L^1$ and $\bbE[ \frac{|X|^p}{X}Y]=0$ hence \eqref{dsuy} implies
\begin{equation}\label{weui}
 \bbE[|X+Y|^p-|X|^p]=\bbE\left[ |X+Y|^p-|X|^p- p \frac{|X|^p}{X}Y \right]\ge c_p\bbE\left[ \min( |X|^{p-2}Y^2, |Y|^p)\right].
\end{equation}
By setting $\delta:= \sqrt{\frac{\bbE[ |Y+X|^p-|X|^p]}{ \bbE[|X|^p] }}\le 1$
we have (using \eqref{weui} in the second line)
\begin{equation}\begin{split}
\bbE[|Y|^p\ind_{|Y|\le  \delta |X|}]&\le \delta^p \bbE[|X|^p] \\
\bbE[|Y|^p\ind_{|Y|> \delta |X|}]&\le \bbE\left[ \delta^{2-p}\min( |X|^{2-p} Y^2, |Y|^p)
\right]\le \frac{\delta^{2-p}}{c_p} \bbE[|X+Y|^p-|X|^p]
\end{split}\end{equation}
and we conclude by adding the two inequalities, obtaining \eqref{dfty} with $C_p=1+\frac{1}{c_p}$.
\end{proof}

\subsection{Proof of Proposition \ref{sideresult}}\label{secside}
We fix $\delta$ small and assume without loss of generality that
\begin{equation}\label{assumep}
 p\in  \left(1+\frac{2}{d}-\frac{\delta}{2d},1+\frac{2}{d} \right).
\end{equation}
We set $u_n:= \beta_n-\beta_c$. We have, for any $k\ge 1$.
\begin{equation}
 \varlimsup_{n\to \infty} \| W^{\beta_c+u_n}_n- W^{\beta_c}_\infty \|_p \le  \varlimsup_{n\to \infty}\| W^{\beta_c+u_n}_n- W^{\beta_c+u_n}_k\|_p +   \varlimsup_{n\to \infty}\| W^{\beta_c+u_n}_k-W^{\beta_c}_k\|_p+ \|W^{\beta_c}_k-W^{\beta_c}_\infty \|_p.
\end{equation}
The second term in the r.h.s.\ is equal to zero by continuity. Hence taking the limit when $k\to \infty$ we obtain that
\begin{equation}
  \varlimsup_{n\to \infty} \| W^{\beta_c+u_n}_n- W^{\beta_c}_\infty \|_p \le
  \varlimsup_{k\to \infty}\varlimsup_{n\to \infty}\| W^{\beta_c+u_n}_n- W^{\beta_c+u_n}_k\|_p
  + \lim_{k\to \infty} \|W^{\beta_c}_k-W^{\beta_c}_\infty \|_p.
\end{equation}
The last term is equal to zero by Theorem \ref{leB}, hence to prove   Proposition \ref{sideresult}
we only need to show that
\begin{equation}\label{ytuio}
 \lim_{k\to \infty}  \varlimsup_{n\to  \infty}\| W^{\beta_c+u_n}_n- W^{\beta_c+u_n}_k\|_p =0.
\end{equation}
As  $W^{\beta_c+u_n}_n- W^{\beta_c+u_n}_k$ is a martingale increment, Lemma \ref{otptilem} applies and we have
\begin{equation}
\| W^{\beta_c+u_n}_n- W^{\beta_c+u_n}_k\|_p \le C^{1/p}_p \|W^{\beta_c+u_n}_k\|_p^{1-\frac{p}{2}} \sqrt{ \bbE\left[ (W^{\beta_c+u_n}_n)^p- (W^{\beta_c+u_n}_k)^p  \right]},
\end{equation}
if the following condition is satisfied
\begin{equation}\label{tyupd}
\bbE\left[ (W^{\beta_c+u_n}_n)^p- (W^{\beta_c+u_n}_k)^p  \right]\le \|W^{\beta_c+u_n}_k\|^p_p.
\end{equation}
Hence to prove \eqref{ytuio} it is sufficient to show that
\begin{equation}\label{treup}
\lim_{n\to \infty} \bbE\left[ \left(W^{\beta_c+u_n}_n\right)^p \right] =\bbE\left[ \left(W^{\beta_c}_\infty\right)^p  \right].
\end{equation}
Indeed \eqref{treup} implies that
\begin{equation}
\lim_{n\to \infty} \bbE\left[ (W^{\beta_c+u_n}_n)^p- (W^{\beta_c+u_n}_k)^p  \right]
= \bbE\left[ (W^{\beta_c}_\infty)^p- (W^{\beta_c}_k)^p  \right]
\end{equation}
and hence \eqref{tyupd} is satisfied if $k$ is large enough and $n\ge n_0(k)$. Furthermore \eqref{treup} also implies that
\begin{equation}
\lim_{n\to \infty} \|W^{\beta_c+u_n}_k\|_p^{1-\frac{p}{2}} \sqrt{ \bbE\left[ (W^{\beta_c+u_n}_n)^p- (W^{\beta_c+u_n}_k)^p  \right]}= \|W^{\beta_c}_k\|_p^{1-\frac{p}{2}} \sqrt{ \bbE\left[ (W^{\beta_c}_\infty)^p- (W^{\beta_c}_k)^p  \right]}
\end{equation}
and the convergence of the r.h.s.\ to zero as $k\to \infty$ follows from Theorem \ref{leB}.
To conclude let us prove \eqref{treup}.
Using Lemma \ref{compapinningz} we have
\begin{equation}\label{lafiin}\begin{split}
 \bbE\left[ \left(W^{\beta_c+u_n}_n\right)^p-\left(W^{\beta_c}_n\right)^p  \right]&\le \bbE\left[ (W^{\beta_c}_\infty)^p\right] \sum_{k\ge 1} \sum_{0<i_1<\dots<i_k\le n} (C u_n^{\frac{p}{2}})^{k} \prod_{j=1}^k\cK^{(p)}(i_{j}-i_{j-1})
 \\ &\le \bbE\left[ (W^{\beta_c}_\infty)^p\right]   \sum_{k\ge 1}  \left(  C u^{\frac p 2 }_n \sum_{\ell =1}^n \mathcal K^{(p)}(\ell) \right)^k
\end{split}\end{equation}
 where the second line is obtained by replacing the condition $i_k\le n$ by a softer one, that is,  all increments in the sequence $(i_j)_{j=1}^k$ are smaller than $n$.
 Proposition \ref{crucialprim} together with the assumptions \eqref{assumep} and $u_n\le n^{-\delta}$ implies that for $n$ sufficiently large we have
 \begin{equation}
 \limsup_{n\to \infty} \log\left( u^{p/2}_n \sum_{\ell =1}^n \mathcal K^{(p)}(\ell)\right) \le d+2-dp-\frac{\delta p}{2}<0
 \end{equation}
and thus in particular
  $\lim_{n\to \infty} u^{\frac p 2 }_n \sum_{\ell =1}^n \mathcal K^{(p)}(\ell)= 0$.
 Using this in \eqref{lafiin} implies that
 $\lim_{n\to \infty} \bbE\left[ \left(W^{\beta_c+u_n}_n\right)^p-\left(W^{\beta_c}_n\right)^p  \right]=0,$
  concluding the proof.
\qed

\medskip

{\bf Ackowledgements:} The author is deeply indebted to S. Junk for numerous enlightening discussions. He wishes to thank the referees for their helpful comments to improve the quality of the manuscript.
He is also grateful to F. Caravenna and R. Sun for helpful comments on an early draft of the manuscript.  He acknowledges the support of a productivity grant from CNPq and a CNE grant from FAPERj.

\bibliographystyle{plain}
\bibliography{ref.bib}

\end{document}